\documentclass[final]{siamltex}

\usepackage{epsfig,amssymb,latexsym}
\usepackage{amsfonts,psfrag,amsmath,bbm,color}
\usepackage{cancel}
\usepackage{siunitx}
\usepackage{comment}
\usepackage{mathrsfs}
\usepackage{graphicx}
\usepackage{textcomp}
\usepackage{multirow}
\usepackage{enumerate}
\usepackage{cancel}
\usepackage{algpseudocode}
\usepackage{caption}  
\usepackage{subcaption}
\usepackage{url}
\usepackage{rotating}
\usepackage{slashbox}
\usepackage{hyperref}
\usepackage{xcolor}
\usepackage{ulem}
\usepackage{mathrsfs}

\usepackage{placeins}

\def\grad{{\nabla}}
\usepackage[ruled,vlined]{algorithm2e}
\usepackage{geometry}
\vbadness=\maxdimen

\newcommand{\footremember}[2]{%
    \footnote{#2}
    \newcounter{#1}
    \setcounter{#1}{\value{footnote}}%
}
\newcommand{\footrecall}[1]{%
    \footnotemark[\value{#1}]%
} 

\usepackage{tikz}
\usetikzlibrary{shapes.geometric, arrows}
\tikzstyle{startstop} = [rectangle, rounded corners, minimum width=1cm, minimum height=1cm,text centered, draw=black]
\tikzstyle{io} = [trapezium, trapezium left angle=70, trapezium right angle=110, minimum width=1cm, minimum height=1cm, text centered, draw=black, fill=blue!30]
\tikzstyle{method} = [rectangle, rounded corners, minimum width=1cm, minimum height =1cm, text centered, draw=black]
\tikzstyle{process} = [rectangle, minimum width=1cm, minimum height=1cm, text centered, draw=black]
\tikzstyle{decision} = [diamond, minimum width=0.5cm, minimum height=0.5cm, text centered, draw=black, fill=green!30]
\tikzstyle{arrow} = [thick,->,>=stealth]

\usepackage{amscd}

\graphicspath{{./figs/}}

\newcommand{\lp}{\left(}
\newcommand{\rp}{\right)}

\newcommand{\lab}{<\hspace{-1mm}}
\newcommand{\rab}{\hspace{-1mm}>}

\newtheorem{remark}{Remark}[section]

\def\PP{{{\rm l}\kern - .15em {\rm P} }}
\def\PN2{{\PP_{N}-\PP_{N-2}}}

\newcommand{\E}{\mathbbm{E}}

\newcommand{\cD}{\mathcal{D}}

\newcommand{\bphi}{\boldsymbol{\varphi}}

\newcommand{\bif}{\textbf{\textit{f}}\hspace{-0.5mm}}
\newcommand{\ba}{\boldsymbol{a}}

\newcommand{\bchi}{\pmb{\chi}}

\newcommand{\be}{\boldsymbol{e}}

\newcommand{\bg}{\textbf{\textit{g}}}

\newcommand{\bH}{\boldsymbol{H}}

\newcommand{\bL}{\boldsymbol{L}}

\newcommand{\hl}{{\hat{l}}}
\newcommand{\hp}{{\hat{p}}}
\newcommand{\hq}{{\hat{q}}}

\newcommand{\bR}{\boldsymbol{R}}

\newcommand{\bu}{\boldsymbol{u}}

\newcommand{\bv}{\boldsymbol{v}}

\newcommand{\bhu}{\hat{\boldsymbol{u}}}
\newcommand{\bhv}{\hat{\boldsymbol{v}}}
\newcommand{\bnh}{\hat{\textbf{\textit{n}}}}
\newcommand{\hlam}{\hat{\lambda}}

\newcommand{\bV}{\boldsymbol{V}}

\newcommand{\bw}{\boldsymbol{w}}
\newcommand{\bW}{\boldsymbol{W}}
\newcommand{\bhw}{\hat{\boldsymbol{w}}}

\newcommand{\btu}{\tilde{\boldsymbol{u}}}

\newcommand{\bx}{\boldsymbol{x}}
\newcommand{\bX}{\boldsymbol{X}}

\newcommand{\bY}{\boldsymbol{Y}}
\newcommand{\by}{\boldsymbol{y}}

\newcommand{\bGamma}{\boldsymbol{\Gamma}}

\newcommand{\red}[1]{{\color{red}#1}}
\newcommand{\blue}[1]{{\color{blue}#1}}

\newcommand{\deleted}[1]{{}}

\begin{document}
\title{Efficient, Accurate, and Robust Penalty-Projection Algorithm for Parameterized Stochastic Navier-Stokes Flow Problems}

\author{
Neethu Suma Raveendran\footremember{uabm}{D\MakeLowercase{epartment of} M\MakeLowercase{athematics}, U\MakeLowercase{niversity of} A\MakeLowercase{labama at} B\MakeLowercase{irmingham}, AL 35294, USA; P\MakeLowercase{artially supported by the} N\MakeLowercase{ational} S\MakeLowercase{cience} F\MakeLowercase{oundation} (NSF) \MakeLowercase{grant} DMS-2425308.}%
\and
Md. Abdul Aziz\footrecall{uabm}
\and Sivaguru S. Ravindran\footremember{uah}{D\MakeLowercase{epartment of} M\MakeLowercase{athematical} S\MakeLowercase{ciences}, U\MakeLowercase{niversity of} A\MakeLowercase{labama in} H\MakeLowercase{untsville}, H\MakeLowercase{untsville}, AL 35899, USA.}\and Muhammad Mohebujjaman\footnote{C\MakeLowercase{orrespondence: mmohebuj@uab.edu}}\hspace{1mm}\footrecall{uabm}
 }

\maketitle

\begin{abstract}
This paper presents and analyzes a fast, robust, efficient, and optimally accurate fully discrete splitting algorithm for the Uncertainty Quantification (UQ) of parameterized Stochastic Navier-Stokes Equations (SNSEs) flow problems those occur in the convection-dominated regimes. The time-stepping algorithm is an implicit backward-Euler linearized method, grad-div and Ensemble Eddy Viscosity (EEV) regularized, and split using discrete Hodge decomposition. Additionally, the scheme's sub-problems are all designed to have different Right-Hand-Side (RHS) vectors but the same system matrix for all realizations at each time-step. The stability of the algorithm is rigorously proven, and it has been shown that appropriately large grad-div stabilization parameters vanish the splitting error. The proposed UQ algorithm is then combined with the Stochastic Collocation Methods (SCMs). Several numerical experiments are given to verify this superior scheme's predicted convergence rates and performance on benchmark problems for high expected Reynolds numbers ($Re$).
\end{abstract}

{\bf Key words.} Finite element method, uncertainty quantification, penalty-projection method, fast ensemble calculation

\medskip
{\bf Mathematics Subject Classifications (2000)}: 65M12, 65M15, 65M60, 76M10

\pagestyle{myheadings}
\thispagestyle{plain}

\markboth{\MakeUppercase{A penalty-projection efficient algorithm for NSE flow ensemble}}{\MakeUppercase{ N. S. Raveendran, M. A. Aziz, S. S. Ravindran, and M. Mohebujjaman}}

\section{Introduction}
  The time-dependent, incompressible, homogeneous Newtonian SNSEs flow is represented by the following dimensionless stochastic non-linear partial differential equations \cite{gunzburger2019evolve}:  \begin{align}
\bu_{t}+\bu\cdot\nabla \bu-\nabla\cdot\left(\nu(\bx,\omega)\nabla \bu\right)+\nabla p &=  \bif(\bx,t,\omega), \hspace{2mm}\text{in}\hspace{2mm} (0,T]\times \cD ,\label{momentum}\\
\nabla\cdot \bu & = 0, \hspace{14mm}\text{in}\hspace{2mm}(0,T]\times \cD, \label{incompressibility}\\
\bu(\bx,t,\omega)&=\bg(\bx,\omega)\hspace{6mm}\text{in}\hspace{2mm}(0,T]\times \partial\cD,\label{bc-condition}\\
\bu(\bx,0,\omega)& = \bu^0(\bx,\omega),\hspace{3mm}\text{in}\hspace{2mm}\cD.\label{nse-initial}
\end{align}
Where, the simulation end time is represented by $T>0$, $\bx$ the spatial variable, $t$ the time variable, $\bif$ the external force, $\bg$ is prescribed for the Dirichlet boundary condition, $\bu^0$ the initial condition, a convex polygon or polyhedron domain $\cD\subset \mathbb{R}^d\ (d=2,3)$, and $\partial\cD$ is the boundary of $\cD$. 
$\nu(\bx,\omega)$ is the random viscosity which is modeled with $\omega\in\Omega$, the set of all outcomes of a complete probability space $(\Omega,\mathcal{F},P)$
where $\mathcal{F}\subset 2^\Omega$ is the $\sigma$-algebra of the events, and $P$ is a probability measure. We want to find the velocity field $\bu:\Lambda\rightarrow\mathbb{R}^d$, and the modified pressure $p:\Lambda\rightarrow\mathbb{R}$, where $\Lambda:=(0,T]\times \cD\times \Omega$. To make our analysis simple, we consider no-slip boundary condition ($\bg=\textbf{0}$).

Models depend on various key parameters that are not possible to measure perfectly, are highly sensitive to the noise, and exhibit the ``butterfly effect" \cite{lorenz1963deterministic} in the prediction. A common treatment for these issues is that the models need to realize a large size ensemble of parameter values. To solve \eqref{weak_formulation_final-1}-\eqref{weak_formulation_final-2}, we consider the following set of $J$ (number of realizations, or the number of stochastic collocation points, and commonly, $J\gg 1$) ensemble computation:\begin{eqnarray}
\bu_{j,t}+\bu_j\cdot\nabla \bu_j-\nabla\cdot\left(\nu_j(\bx) \nabla \bu_j\right)+\nabla p_j &= & \bif_j(\bx,t), \hspace{2mm}\text{in}\hspace{2mm}(0,T]\times \cD, \label{gov1}\\
\nabla\cdot \bu_j & =& 0, \hspace{11mm}\text{in}\hspace{2mm}(0,T]\times \cD,\label{gov2}
\end{eqnarray}
where $\bu_j$, and $p_j$, denote the velocity and pressure variables, respectively, for each $j=1,2,\cdots,J$, for different viscosity $\nu_j$, and/or body forces $\bif_j$,  and/or the initial conditions $\bu_j^0$, and/or the $J$ different boundary conditions. We specified the no-slip boundary conditions for every realization in order to keep our analysis simple. It is assumed that $\nu_j(\bx)\in L^\infty(\cD)$, and $\nu_j(\bx)\ge\nu_{j,\min}>0$, where $\nu_{j,\min}=\min\limits_{\bx\in\cD}\nu_j(\bx)$.


For a reliable and high-fidelity solution in a wide range of applications, such as hydrology \cite{GG11}, spectral methods \cite{LK10}, weather forecasting \cite{L05,LP08}, porous media flow \cite{jiang2021artificial}, magnetohydrodynamics (MHD) \cite{jiang2018efficient,MR17}, surface data assimilation \cite{fujita2007surface}, and sensitivity analyses \cite{MX06}, the computation of an ensemble average solution is frequently utilized. 

The standard ensemble average solution is computed by running simulations for all realizations independently and then taking the average of them, e.g., the Monte Carlo (MC) methods. Thus, the total cost of solving \eqref{gov1}-\eqref{gov2} $=J\times$ One Navier-Stokes Equations (NSEs) simulation. 

To reduce this huge computational cost by an order of magnitude, a breakthrough idea was proposed by Jiang and Layton in 2014 \cite{JL14}. It designs the scheme so that at each time-step the system matrix of the saddle-point problems remains exactly the same for all sample points but with different RHS vectors. Thus, at each time-step, it requires solving the following block saddle-point linear system for the velocity and pressure nodal values\begin{align}
\begin{pmatrix}\mathbb{A} & \mathbb{B}^T\\\mathbb{B} & \mathcal{O}\end{pmatrix}\bigg(\begin{matrix}
\textbf{U}_1\\
     \textbf{P}_1
\end{matrix}
\bigg\rvert\begin{matrix}
\textbf{U}_2\\
\textbf{P}_2
\end{matrix}\bigg\rvert\begin{matrix}
    \cdots\\
    \cdots
\end{matrix}\bigg\rvert\begin{matrix}
\textbf{U}_J\\
\textbf{P}_J
\end{matrix}\bigg)=\bigg(\begin{matrix}
    \textbf{F}_1\\\textbf{0}
\end{matrix}\bigg\rvert\begin{matrix}
    \textbf{F}_2\\\textbf{0}
\end{matrix}\bigg\rvert\begin{matrix}
    \cdots\\
    \cdots
\end{matrix}\bigg\rvert\begin{matrix}
    \textbf{F}_J\\\textbf{0}
\end{matrix}\bigg),\label{sparse-system}
 \end{align} where the block coefficient matrix $\begin{pmatrix}\mathbb{A} & \mathbb{B}^T\\\mathbb{B} & \mathcal{O}\end{pmatrix}$ is independent of the index $j$, $\mathbb{A}$ is the block-matrix corresponding to the velocity nodal values, $\mathbb{B}$ represents the gradient operator, $\mathbb{B}^T$ its adjoint, and $\mathcal{O}$ represents a zero block-matrix. Also, $\textbf{U}_j$, and $\textbf{P}_j$ are the nodal vectors of velocity, and pressure fields, respectively, and $\textbf{F}_j$ is the appropriate RHS vector for $j=1,2,\cdots,J$. 
 
 Though it saves a huge amount of computational cost for the UQ of SNSEs flow problems, a fundamental hurdle in solving the saddle problem is that the velocity and pressure are coupled. A popular remedy for this is to use splitting methods, which were introduced first by Chorin and Temam \cite{chorin1968numerical, temam1969approximation}. But the accuracy reduces due to the spitting errors, which are caused by several issues. For example, in the `Projection Methods,' two velocity solutions are computed where one does not satisfy the incompressibility condition, and the other does not satisfy the boundary condition and is often considered a non-physical pressure boundary condition. A Penalty Projection (PP) method is analyzed and tested in \cite{linke2017connection} for deterministic NSEs that show that for an appropriately large penalty parameter, the splitting error vanishes. This PP method is being used in an efficient and robust UQ algorithm for the Stochastic MHD convection-dominated flow problems in \cite{mohebujjaman2024efficient}. 
 
 Since the NSEs are the basis for simulating flows of computational modeling and other multi-physics flow problems, it is important to propose, analyze, and test a robust, fast, and efficient splitting algorithm for the parameterized SNSEs without compromising accuracy. To this end, we propose a robust, accurate, and efficient grad-div Stabilized PP with EEV (SPP-EEV) regularized method for SNSEs, which incorporates the breakthrough idea presented in \cite{JL14}  to solve the \eqref{gov1}-\eqref{gov2} ensemble system.
Throughout this paper, a uniform time-step size $\Delta t$ is considered. We define $t^n=n\Delta t$ for $n=0, 1, \cdots$, then compute the $J$ solutions as below: For $j=1,\cdots,J$,\\
Step 1: Find $\bhu_{j,h}^{n+1}$:
\begin{align}
\frac{\bhu_{j,h}^{n+1}}{\Delta t}+<\bhu_h>^n\cdot\nabla \bhu_{j,h}^{n+1}-\nabla\cdot\left(\Bar{\nu}\nabla \bhu_{j,h}^{n+1}\right)&-\gamma\nabla(\nabla\cdot\bhu_{j,h}^{n+1})-\nabla\cdot\left( 2\nu_T(\bhu^{'}_{h},t^n)\nabla \bhu_{j,h}^{n+1}\right)\nonumber\\&= \bif_{j}(t^{n+1})+\frac{\tilde{\bu}_{j,h}^n}{\Delta t}-\bhu_{j,h}^{'n}\cdot\nabla \bhu_{j,h}^n+\nabla\cdot\left(\nu_j^{'}\nabla \bhu_{j,h}^{n}\right),\label{scheme1}\\\bhu_{j,h}^{n+1}\big|_{\partial\cD}&=0.\label{incom1}
\end{align}
Step 2: Find $\tilde{\bu}_{j,h}^{n+1}$, and $\hp_{j,h}^{n+1}$:
\begin{align}
    \frac{\tilde{\bu}_{j,h}^{n+1}}{\Delta t}+\nabla \hp_{j,h}^{n+1}&=\frac{\bhu^{n+1}_{j,h}}{\Delta t},\label{step-2-first}\\
    \nabla\cdot\tilde{\bu}_{j,h}^{n+1}&=0,\label{step-2-second}\\
\tilde{u}_{j,h}^{n+1}\cdot\bnh\big|_{\partial\cD}&=0.\label{normal-eqn}
\end{align}
Where $\bhu_{j,h}^n$, and $\hp_{j,h}^{n}$ denote approximations of $\bu_j(\cdot,t^n)$, and $p_j(\cdot,t^n)$, respectively, $\tilde{\bu}_{j,h}^{n}$ is the projection of $\bhu_{j,h}^n$ onto the divergence-free space, the grad-div parameter $\gamma>0$, and the outward unit normal vector $\bnh$. The ensemble mean and fluctuation are defined as follows:
\begin{align*}
<\bhu_h>^n:=\frac{1}{J}\sum\limits_{j=1}^{J}\bhu_{j,h}^n, \bhu_{j,h}^{'n}:=\bhu_{j,h}^n-<\bhu_h>^n, \Bar{\nu}:=\frac{1}{J}\sum\limits_{j=1}^{J}\nu_j,\;\text{and}\;\nu_j^{'}:=\nu_j-\Bar{\nu}.   
\end{align*}
Using mixing length, the EEV is defined as below \cite{jiang2015numerical}:
\begin{align} \nu_T(\bhu_h^{'},t^n):=\mu\Delta t(\hl^{n})^2,\hspace{0.5mm}\text{ where}\hspace{1mm}(\hl^{n})^2=\sum_{j=1}^J|\bhu_{j,h}^{'n}|^2 \label{eddy-viscosity}
\end{align} is a scalar quantity, $\mu$ is a calibration constant or tuning parameter, and $|\cdot|$ denotes Euclidean norm of a vector. The EEV is a $O(\Delta t)$ term.  For convection-dominated flow, the EEV term helps the model in reducing the numerical instability that arises for under resolved meshes. The concept of EEV is adopted from turbulence modeling and is used widely  \cite{jiang2015higher,jiang2015analysis,jiang2015numerical,Mohebujjaman2022High, mohebujjaman2024efficient, MR17,mohebujjaman2022efficient}. However, the analysis of the algorithms equipped with EEV is scarce.  

Note that in Step 1, the no-slip boundary condition is assigned whereas in Step 2, no penetration boundary condition is considered. At each time-step, each sub-problem (in Steps 1-2) has the same coefficient matrix for all realizations. The global system matrix assembly, factorization (if a direct solution is employed), and preconditioner building (especially for the iterative solver) are thus required only once per time-step, saving enormous amounts of computer memory. Additionally, the algorithm can benefit from block linear solvers. The idea of designing such an ensemble algorithm has been implemented to heat equation \cite{Unconditionally2021Fiordilino}, NSEs \cite{jiang2015higher,jiang2017second,jiang2015numerical,neda2016ensemble}, parameterized flow problems \cite{GJW18, GJW17}, MHD \cite{jiang2018efficient,MR17,Mohebujjaman2022High,mohebujjaman2022efficient,mohebujjaman2024efficient},  and turbulence modeling \cite{JKL15}. 

In Step 1, we need to solve a much smaller system compare to \eqref{sparse-system}. On the other hand, in Step 2, however, we need to solve a similar $2\times 2$ block system as in \eqref{sparse-system}, it is much easir to solve, since the non-linear and the viscous terms are absent and the block-matrix is symmetric positive definite.


Thus, using a finite element spatial discretization, in this paper, we investigate the new efficient, accurate and robust SPP-EEV ensemble method  in a fully discrete form. The stability of the SPP-EEV Algorithm is proven and shown that as $\gamma\rightarrow\infty$, it's solution converges to the solution of an equivalent Coupled EEV (Coupled-EEV) regularized  method for \eqref{gov1}-\eqref{gov2}. Finally, we connect the SPP-EEV and Coupled-EEV schemes to the sparse-grid SCMs \cite{FTW2008} which require fewer realizations compare to the MC methods.  

An ensemble based penalty-projection method is studied in \cite{yuan2020ensemble} without EEV which is not long-time stable and blows-up for more realistic convection-dominated flows, see FIG. 9 in \cite{jiang2015numerical}, see Fig. 6 in \cite{mohebujjaman2024efficient}, and in this paper see Fig. \ref{RLDC-energy:mu-varies}. To our knowledge, the SPP-EEV is a novel parameterized algorithm for fast, robust, efficient, and optimally accurate computation of UQ of SNSEs flow problems. Moreover, it can be an enabling tool for robust computation in fluid dynamics, electromagnetism, and plasma dynamics.

The report's remaining sections are arranged as follows: We include the necessary notations and mathematical preliminaries in Section \ref{notation-prelims} to ensure a smooth analysis. In Section \ref{fully-discrete-scheme}, we present the Coupled-EEV algorithm, and state it's stability and convergent theorems. We then propose the novel more efficient SPP-EEV algorithm in fully-discrete setting in Section \ref{penalty-projection}, where we prove the stability and convergence theorems. It is shown that as $\gamma\rightarrow\infty$, the SPP-EEV scheme converges to the Coupled-EEV scheme. In Section \ref{scm}, we provide a brief introduction of the SCMs, combine it with the SPP-EEV and Coupled-EEV methods for solving SNSEs further efficiently. Several numerical experiments are given in Section \ref{numerical-experiment}, which support the theory.

\section{Notation and preliminaries}\label{notation-prelims}
 The usual $L^2(\cD)$ norm and inner product are denoted by $\|.\|$ and $(\cdot,\cdot)$, respectively. Similarly, the $L^p(\cD)$ norms and the Sobolev $W_p^k(\cD)$ norms are $\|.\|_{L^p}$ and $\|.\|_{W_p^k}$, respectively for $k\in\mathbb{N},\hspace{1mm}1\le p\le \infty$. Sobolev space $W_2^k(\cD)$ is represented by $H^k(\cD)$ with norm $\|.\|_k$. The vector-valued spaces are $$\bL^p(\cD)=(L^p(\cD))^d, \hspace{1mm}\text{and}\hspace{1mm}\bH^k(\cD)=(H^k(\cD))^d.$$
For $\bX$ being a normed function space in $\cD$, $L^p(0,T;\bX)$ is the space of all functions defined on $(0,T]\times\cD$ for which the following norm
\begin{align*}
\|\bu\|_{L^p(0,T;\bX)}=\lp\int_0^T\|\bu\|_{\bX}^pdt\rp^\frac{1}{p}\hspace{-2mm},\hspace{2mm}p\in[1,\infty)
\end{align*}
is finite. For $p=\infty$, the usual modification is used in the definition of this space. The natural function spaces for our problem are
\begin{align*}
    \bX:&=\bH_0^1(\cD)=\{\bv\in \bL^2(\cD) :\nabla \bv\in L^2(\cD)^{d\times d}, \bv=0 \hspace{2mm} \mbox{on}\hspace{2mm}   \partial \cD\},\\
    \bY:&=\{\bv\in  L^2(\cD):\nabla\cdot\bv\in L^2(\cD),\bv\cdot\bnh\big|_{\partial\cD}=0\},\\
    Q:&=L_0^2(\cD)=\{ q\in L^2(\cD): \int_\cD q\hspace{1mm}d\bx=0\}.
\end{align*}
Recall the Poincar\'e inequality holds in $\bX$: There exists $C$ depending only on $\cD$ satisfying for all $\bphi\in \bX$,
\[
\| \bphi \| \le C \| \nabla \bphi \|.
\]
The divergence free velocity space is given by
$$\bV:=\{\bv\in \bX:(\nabla\cdot \bv, q)=0, \forall q\in Q\}.$$
We define the skew symmetric trilinear form $b^*:\bX\times \bX\times \bX\rightarrow \mathbb{R}$ by
	\[
	b^*(\bu,\bv,\bw):=\frac12(\bu\cdot\nabla \bv,\bw)-\frac12(\bu\cdot\nabla \bw,\bv),\hspace{3mm} 
	\]
 By the divergence theorem \cite{jiang2015higher}, it can be shown \begin{align}
     b^*(\bu,\bv,\bw)=(\bu\cdot\nabla \bv,\bw)+\frac12(\nabla\cdot\bu,\bv\cdot\bw),\;\;\;\forall\bu,\bv,\bw\in\bX.\label{trilinear-identitiy}
 \end{align}Recall from \cite{L08, lee2011error, linke2017connection} that for any $\bu,\bv,\bw\in 
		\bX$
	\begin{align}
		b^*(\bu,\bv,\bw)&\leq C(\cD)\|\nabla \bu\|\|\nabla \bv\|\|\nabla \bw\|,\label{nonlinearbound}
	\end{align}	
 and additionally, if $\bv\in \bL^\infty(\Omega)$, and $\nabla\bv\in\bL^3(\Omega)$, then 
 \begin{align}
    b^*(\bu,\bv,\bw)\leq C(\cD)\|\bu\|\left(\|\nabla\bv\|_{L^3}+\|\bv\|_{L^\infty}\right)\|\nabla\bw\|. \label{nonlinearbound3}
 \end{align}The following inequality is proved in Appendix \ref{proof-basic-eqn}, will be used frequently
\begin{align}
\|\bu\cdot\nabla\bv\|&\le\||\bu|\nabla\bv\|.\label{basic-ineq}
\end{align}
The conforming finite element spaces are denoted by $\bX_h\subset \bX$ and  $Q_h\subset Q$, and we assume a regular triangulation $\tau_h(\cD)$, where $h$ is the maximum triangle diameter.   We assume that $(\bX_h,Q_h)$ satisfies the usual discrete inf-sup condition
\begin{eqnarray}
\inf_{q_h\in Q_h}\sup_{\bv_h\in \bX_h}\frac{(q_h,\grad\cdot \bv_h)}{\|q_h\|\|\grad \bv_h\|}\geq\beta>0,\label{infsup}
\end{eqnarray}
where $\beta$ is independent of $h$. We assume that there exists a finite element space  $\bY_h\subset\bY$. The space of discretely divergence free functions is defined as 
\begin{align*}
    \bV_h:=\{\bv_h\in \bX_h:(\nabla\cdot \bv_h,q_h)=0,\hspace{2mm}\forall q_h\in Q_h\}.
\end{align*}


We will assume the mesh is sufficiently regular for the inverse inequality to hold. The following lemma for the discrete Gr\"onwall inequality was given in \cite{HR90}.
\begin{lemma}\label{dgl}
		Let $\Delta t$, $\mathcal{E}$, $a_n$, $b_n$, $c_n$, $d_n$ be non-negative numbers for $n=1,\cdots, M$ such that
		$$a_M+\Delta t \sum_{n=1}^Mb_n\leq \Delta t\sum_{n=1}^{M-1}{d_na_n}+\Delta 
		t\sum_{n=1}^Mc_n+\mathcal{E}\hspace{3mm}\mbox{for}\hspace{2mm}M\in\mathbb{N},$$
		then for all $\Delta t> 0,$
		$$a_M+\Delta t\sum_{n=1}^Mb_n\leq \mbox{exp}\left(\Delta t\sum_{n=1}^{M-1}d_n\right)\lp\Delta 
		t\sum_{n=1}^Mc_n+\mathcal{E}\rp\hspace{2mm}\mbox{for}\hspace{2mm}M\in\mathbb{N}.$$
	\end{lemma}
To simplify the notation, denote $\alpha_{\min}:=\min\limits_{1\le j\le J}\alpha_j$, where $\alpha_j:=\Bar{\nu}_{\min}-\|\nu_j^{'}\|_\infty$, for $j=1,2,\cdots\hspace{-0.35mm}, J$, and $\Bar{\nu}_{\min}:=\min\limits_{\bx\in\cD}\Bar{\nu}(\bx)$.  We assume that the data has no outlier and that observations are close enough to the mean, so $\alpha_j>0$ holds.
 \section{Fast, Robust, and Efficient Coupled-EEV algorithm for SNSEs}\label{fully-discrete-scheme}
In this section, we present a velocity and pressure coupled, fully discrete, efficient, linear extrapolated backward-Euler method in Algorithm \ref{coupled-alg} which is EEV and grad-div regularized finite element time-stepping algorithm for the parameterized SNSEs. This Algorithm \ref{coupled-alg} is a variation of Algorithm 1 in \cite{raveendran2024two} (without grad-div term), and also variation of an algorithm given in (2.1) in \cite{gunzburger2019efficient} (without grad-div and EEV terms). \\ 
\begin{algorithm}[H]\label{coupled-alg}
  \caption{Coupled-EEV scheme} Given $T>0$, $\gamma>0$, $\Delta t>0$, $\bu_j^0\in  \bV_h$, and $\bif_{j}\in$ $ L^\infty\left( 0,T;\bH^{-1}(\cD)\right)$, set $M=T/\Delta t$. If $\alpha_j>C\|\nabla\cdot\bu_{j,h}^{'n}\|_{L^\infty}$, choose  $\mu\ge \frac{\alpha_j}{\Delta t\left(C\alpha_j^2-\|\nabla\cdot\bu_{j,h}^{'n}\|_{L^\infty}^2\right)}$, $\forall j=1,2,\cdots\hspace{-0.35mm},J$. Then, for $n=1,\cdots\hspace{-0.35mm},M-1$, compute: Find $(\bu_{j,h}^{n+1}, p_{j,h}^{n+1})\in \bX_h\times Q_h$ satisfying, $\forall$ $(\bchi_{h},q_{h})\in \bX_h\times Q_h$:
 \begin{align}
&\Big(\frac{\bu_{j,h}^{n+1}-\bu_{j,h}^n}{\Delta t},\bchi_{h}\Big)+b^*\big(\hspace{-1mm}<\bu_h>^n, \bu_{j,h}^{n+1},\bchi_{h}\big)+\big(\Bar{\nu}\nabla \bu_{j,h}^{n+1},\nabla\bchi_{h}\big)+(\gamma\nabla\cdot\bu_{j,h}^{n+1}-p_{j,h}^{n+1},\nabla\cdot\bchi_{h})\nonumber\\&+\left( 2\nu_T(\bu^{'}_{h},t^n)\nabla \bu_{j,h}^{n+1},\nabla\bchi_h\right)= \big(\bif_{j}(t^{n+1}),\bchi_{h}\big)-b^*(\bu_{j,h}^{'n}, \bu_{j,h}^n,\bchi_{h})-\big(\nu_j^{'}\nabla \bu_{j,h}^{n},\nabla\bchi_{h}\big),
\label{couple-eqn-1}\\\nonumber\\&\big(\nabla\cdot\bu_{j,h}^{n+1},q_{h}\big)=0,\label{couple-incompressibility}
\end{align}
\end{algorithm}
\noindent where $\bu_{j,h}^n$, and $p_{j,h}^{n}$ denote approximations of $\bu_j(\cdot,t^n)$, and $p_j(\cdot,t^n)$, respectively, and EEV is defined as \begin{align} \nu_T(\bu_h^{'},t^n):=\mu\Delta t(l^{n})^2,\hspace{0.5mm}\text{ where}\hspace{1mm}(l^{n})^2=\sum_{j=1}^J|\bu_{j,h}^{'n}|^2. \label{eddy-viscosity2}
\end{align} 

The Algorithm \ref{coupled-alg} is efficient because it is presented in a clever way so that at each time-step, for all the realizations, we need to solve a system as in \eqref{sparse-system}, where the block-matrix $\mathbb{A}$ is for the mass, non-linear, and stiffness matrices. Therefore, it reduces huge computational complexity in terms of computing time and computer memory requirements. Moreover, appropriately large grad-div parameter values will allow the scheme to use \textit{inf-sup} stable weakly divergence-free element, e.g., Taylor-Hood (TH) element \cite{wieners2003taylor}, without unreasonably growing divergence error \cite{case2011connection}. For a given mesh, a pointwise divergence-free element, e.g., Scott-Vogelius (SV) \cite{scott1985norm}, requires more pressure degrees of freedom (dof) than the TH element. The presence of the EEV in the model will provide a long-time stable solution for the convection-dominated flows \cite{mohebujjaman2024efficient}. 
\begin{theorem}
   (Unconditional Stability) Suppose $\bif_{\hspace{0.5mm}j}\in L^\infty(0,T;\bH^{-1}(\cD))$, and $\bu_{j,h}^0\in\bV_h$, then the solutions of Algorithm \ref{coupled-alg} are stable: For all $\Delta t>0$, if $\alpha_j>C\|\nabla\cdot\bu_{j,h}^{'n}\|_{L^\infty}$, choose $\mu\ge \frac{\alpha_j}{\Delta t\left(C\alpha_j^2-\|\nabla\cdot\bu_{j,h}^{'n}\|_{L^\infty}^2\right)}$ and $\gamma>0$, we then have the following stability bound: 
\begin{align}
\|\bu_{j,h}^{M}\|^2+\Bar{\nu}_{\min}\Delta t\|\nabla \bu_{j,h}^{M}\|^2+2 \gamma\Delta t\sum_{n=1}^{M}\|\nabla\cdot\bu_{j,h}^{n}\|^2\le C(data).\label{stability-couple-alg}
\end{align}
\end{theorem}
\begin{proof}
    See the Appendix A in \cite{raveendran2024two}
\end{proof}
\begin{remark}
    In \cite{raveendran2024two}, there is a time-step restriction, which has been avoided without forming $\|\bu_{j,h}^{n+1}-\bu_{j,h}^{n}\|$ in $b^*(\bu_{j,h}^{'n}, \bu_{j,h}^n,\bu_{j,h}^{n+1})$.
\end{remark}
\begin{theorem}(Convergence) Suppose $(\bu_j,p_j)$ satisfying \eqref{momentum}-\eqref{nse-initial} and the following regularity assumptions for $m=\max\{3,k+1\}$
\begin{align*}
   &\bu_j\in L^2(0,T;\bH^{m}(\cD))\cap L^\infty(0,T; \bH^{m}(\cD)),\\
   &\bu_{j,t}\in L^\infty(0,T;\bH^2(\cD))\cap L^2(0,T;\bH^1(\cD)),\bu_{j,tt}\in L^2(0,T;\bL^2(\cD)),
\end{align*}
with $k\ge 2$, then the ensemble solution of the Algorithm \ref{coupled-alg} converges to the true ensemble solution: For $\mu\ge\frac12$, if $\Delta t<\frac{C\alpha_j}{\|\nabla\cdot\bu_{j,h}^{'n}\|^2_{L^\infty}}$ then, the following holds
\begin{align}
    \|<\bu>(T)-<\bu_h>^M\|^2+\alpha_{\min}\Delta t\sum_{n=1}^M\Big\|\nabla\Big(\hspace{-1.1mm}<\bu>(t^n)-<\bu_h>^n\hspace{-1.1mm}\Big)\Big\|^2\nonumber\\\le C\Big(h^{2k}+\Delta t^2+h^{2-d}\Delta t^2+h^{2k-1}\Delta t\Big).\label{convergence-error}
\end{align}
\end{theorem}

\begin{proof}
    Please see the Appendix B in \cite{raveendran2024two}
\end{proof}

\begin{remark}
    It is straightforward to prove the above stability and convergence theorems thus, we omit their proofs. The error estimate in \eqref{convergence-error} shows that in 2D, an optimal convergence can be achieved, however, in 3D, the convergence is sub-optimal. This is because of the EEV term in the model.
\end{remark} 

\begin{lemma}\label{lemma1}
If $\bu_j\in L^\infty(0,T;\bH^2(\cD))$ then for sufficiently small $h$ and $\Delta t$, $\exists$ $C_*\in\mathbb{R}^+$ (which does not depend on $h$, $\Delta t$, and $\gamma$) such that
\begin{align*}
    \max_{1\le n\le M}\Big(\|\nabla \bu_{j,h}^n\|_{L^3}+\|\bu_{j,h}^n\|_{L^\infty}\Big)&\le C_*,\hspace{2mm}\text{for all}\hspace{2mm}j=1,2,\cdots,J.
\end{align*}
\end{lemma}
\begin{proof}
    The proof is given in Appendix \ref{appendix-A}.
\end{proof}

\section{Fast, Robust, and Efficient SPP-EEV algorithm for SNSEs}
\label{penalty-projection}
Now, we present and analyze a more efficient, decoupled, fully discrete, linearly extrapolated backward-Euler method in Algorithm \ref{SPP-FEM}, which is also an EEV and grad-div stabilized finite element time-stepping scheme for computing parameterized SNSEs flow ensemble.

\begin{algorithm}[H]\label{SPP-FEM}
  \caption{SPP-EEV scheme} Given time-step $\Delta t>0$, end time $T>0$, initial conditions $\bhu_j^0=\tilde{\bu}_j^0\in \bY_h\cap \bH^2(\cD)$ and $\bif_{j}\in$ $ L^\infty\left( 0,T;\bH^{-1}(\cD)\right)$ for $j=1,2,\cdots\hspace{-0.35mm},J$. Set $M=T/\Delta t$ and for $n=1,\cdots\hspace{-0.35mm},M-1$, compute:\\
 Step 1: Find $\bhu_{j,h}^{n+1}\in \bX_h$ satisfying for all $\bchi_{h}\in \bX_h$:
 \begin{align}
&\Big(\frac{\bhu_{j,h}^{n+1}-\tilde{\bu}_{j,h}^n}{\Delta t},\bchi_{h}\Big)+b^*\big(\hspace{-1mm}<\hspace{-1mm}\bhu_h\hspace{-1mm}>^n, \bhu_{j,h}^{n+1},\bchi_{h}\big)+\big(\Bar{\nu}\nabla \bhu_{j,h}^{n+1},\nabla\bchi_{h}\big)+\gamma\big(\nabla\cdot\bhu_{j,h}^{n+1},\nabla\cdot \bchi_{h}\big)\nonumber\\&+\Big(2\nu_T(\hat{\bu}^{'}_{h},t^n)\nabla \bhu_{j,h}^{n+1},\nabla \bchi_{h}\Big)= \big(\bif_{j}(t^{n+1}),\bchi_{h}\big)-b^*\big(\bhu_{j,h}^{'n}, \bhu_{j,h}^n,\bchi_{h}\big)-\big(\nu_j^{'}\nabla \bhu_{j,h}^{n},\nabla\bchi_{h}\big).\label{spp-step-1}
\end{align}
Step 2: Find $\left(\tilde{\bu}_{j,h}^{n+1},\hp_{j,h}^{n+1}\right)\in \bY_h\times Q_h$ satisfying for all $(\bv_{h},q_{h})\in \bY_h\times Q_h$:
\begin{align}
\Big(\frac{\tilde{\bu}_{j,h}^{n+1}-\bhu_{j,h}^{n+1}}{\Delta t},\bv_{h}\Big)-\big( \hp_{j,h}^{n+1},\nabla\cdot\bv_{h}\big)&=0,\label{spp-step-2-1}\\
    \big(\nabla\cdot\tilde{\bu}_{j,h}^{n+1},q_{h}\big)&=0.\label{spp-step-2-2}
\end{align}
\end{algorithm}
In Step 1, we solve a $1\times 1$ block system of the form 
\begin{align}
\hat{\mathbb{A}}\big(\begin{matrix}
\hat{\bf{U}}_1
\end{matrix}
\big\rvert\begin{matrix}
	\hat{\bf{U}}_2
\end{matrix}\big\rvert\begin{matrix}
	\cdots
\end{matrix}\big\rvert\begin{matrix}
	\hat{\bf{U}}_J
\end{matrix}\big)=\big(\begin{matrix}
	\hat{\bf{F}}_1
\end{matrix}\big\rvert\begin{matrix}
	\hat{\bf{F}}_2
\end{matrix}\big\rvert\begin{matrix}
	\cdots
\end{matrix}\big\rvert\begin{matrix}
	\hat{\bf{F}}_J
\end{matrix}\big),\label{step-1-bloc}
\end{align}
where $\hat{\mathbb{A}}$ is the coefficient matrix which is independent of the index $j$, and $\hat{\bf{U}}_j$ is the nodal vector for $\bhu_{j,h}^{n+1}$. The size of $\hat{\mathbb{A}}$ is much smaller than that in \eqref{sparse-system}. The Darcy system in Step 2 is essentially the $L^2$ projection of $\btu_{j,h}^{n+1}$ onto the following affine variety of the divergence-free space
$$\{\bv_h\in\bY_h:\nabla\cdot\bv_h=0\}.$$ Here, we solve a $2\times 2$ block system of the form
\begin{align}
\begin{pmatrix}\tilde{\mathbb{A}} & \tilde{\mathbb{B}}^T\\\tilde{\mathbb{B}}& \mathcal{O}\end{pmatrix}\bigg(\begin{matrix}
\tilde{\bf{U}}_1\\
     \hat{\bf{P}}_1
\end{matrix}
\bigg\rvert\begin{matrix}
\tilde{\bf{U}}_2\\
\hat{\bf{P}}_2
\end{matrix}\bigg\rvert\begin{matrix}
    \cdots\\
    \cdots
\end{matrix}\bigg\rvert\begin{matrix}
\tilde{\bf{U}}_J\\
\hat{\bf{P}}_J
\end{matrix}\bigg)=\bigg(\begin{matrix}
    \tilde{\bf{F}}_1\\\textbf{0}
\end{matrix}\bigg\rvert\begin{matrix}
    \tilde{\bf{F}}_2\\\textbf{0}
\end{matrix}\bigg\rvert\begin{matrix}
    \cdots\\
    \cdots
\end{matrix}\bigg\rvert\begin{matrix}
    \tilde{\bf{F}}_J\\\textbf{0}
\end{matrix}\bigg),\label{sparse-system2}
 \end{align}
where $\tilde{\mathbb{A}}$ is for the mass matrix only, $\tilde{\mathbb{B}}$ is the gradient operator, and $\tilde{\mathbb{B}}^T$ is the adjoint of $\tilde{\mathbb{B}}$. Also, $\tilde{\bf{U}}_j$, and $\hat{\bf{P}}_j$ are the nodal vectors of $\tilde{\bu}_{j,h}^{n+1}$, and $\hp_{j,h}^{n+1}$, respectively. The $2\times 2$ block matrix in \eqref{sparse-system2} is symmetric positive definite and much easier to solve the system \eqref{sparse-system2} compared to that in \eqref{sparse-system}. Moreover, taking divergence, and dot product with the normal vector $\bnh$ on both sides of \eqref{step-2-first}, then using the incompressibility constraint, and boundary condition, respectively, to get
\begin{align}
    \Delta \hp_{j,h}^{n+1}&=\frac{1}{\Delta t}\nabla\cdot\bhu^{n+1}_{j,h},\label{step-2-reduced}\\
    \bnh\cdot\nabla \hp_{j,h}^{n+1}\big|_{\partial\cD}&=0.\label{step-2-reduced-bc}
\end{align}
For each realization, the cost of solving \eqref{step-2-reduced}-\eqref{step-2-reduced-bc} is equivalent to solving a Poisson equation at each time-step. Further, we can use a block solve as below:
\begin{align}
\hat{\mathbb{S}}\big(\begin{matrix}
\hat{\bf{P}}_1
\end{matrix}
\big\rvert\begin{matrix}
	\hat{\bf{P}}_2
\end{matrix}\big\rvert\begin{matrix}
	\cdots
\end{matrix}\big\rvert\begin{matrix}
	\hat{\bf{P}}_J
\end{matrix}\big)=\big(\begin{matrix}
	\hat{\bf{G}}_1
\end{matrix}\big\rvert\begin{matrix}
	\hat{\bf{G}}_2
\end{matrix}\big\rvert\begin{matrix}
	\cdots
\end{matrix}\big\rvert\begin{matrix}
	\hat{\bf{G}}_J
\end{matrix}\big),\label{laplace-block}
\end{align}
where $\hat{\mathbb{S}}$ is the stiffness matrix. And, $\hat{\bf{F}}_j, \tilde{\bf{F}}_j$  and $\hat{\bf{G}}_j$, $j=1,2,\cdots,J$ are the appropriate RHS vectors  in \eqref{step-1-bloc}, \eqref{sparse-system2}, and \eqref{laplace-block}, respectively. That is, the block saddle-point system \eqref{sparse-system2} solving cost reduces to the cost of solving block Poisson system \eqref{laplace-block}. Furthermore, $\hat{\mathbb{S}}$ is independent of time, thus, once we construct it for one time-step, can reuse it for all time-steps. Moreover, for a large $\gamma$, the splitting error diminishes.

 Hence, at each time-step, Coupled-EEV, and SPP-EEV approximate the system \eqref{gov1}-\eqref{gov2} solving a single block saddle-point system in \eqref{sparse-system}, and two easier block systems in \eqref{step-1-bloc} and \eqref{laplace-block}, respectively. Therefore, though SPP-EEV requires two block solves instead of one, the computational cost of SPP-EEV method will be less than that of Coupled-EEV method. These ensure the superiority of the SPP-EEV scheme over the Coupled-EEV.

\subsection{Stability Analysis}\label{stability-analysis} This section states and proves the stability and well-posedness of the Algorithm \ref{SPP-FEM}.

\begin{lemma}(Stability)\label{stability-projection}
Let $\big(\bhu_{j,h}^{n+1},\hp_{j,h}^{n+1}\big)$ be the solution of Algorithm \ref{SPP-FEM} and $\textbf{f}_{j}\in L^\infty\left(0,T;\bH^{-1}(\cD)\right)$, and $\bhu_{j,h}^0\in \bV_h$. Then $\forall$ $\Delta t>0$, if $\alpha_j>C\|\nabla\cdot\bhu_{j,h}^{'n}\|_{L^\infty}$, choose $\mu\ge \frac{\alpha_j}{\Delta t\left(C\alpha_j^2-\|\nabla\cdot\bhu_{j,h}^{'n}\|_{L^\infty}^2\right)}$, we then have
\begin{align}
\|\bhu_{j,h}^{M}\|^2+\Bar{\nu}_{\min}\Delta t\|\nabla\bhu_{j,h}^M\|^2+2\gamma\Delta t\sum_{n=1}^{M}\|\nabla\cdot\bhu_{j,h}^{n}\|^2\nonumber\\\le\|\bhu_{j,h}^0\|^2+\Bar{\nu}_{\min}\Delta t\|\nabla \bhu_{j,h}^{0}\|^2+\frac{2\Delta t}{\alpha_j}\sum_{n=1}^{M}\|\bif_{j}(t^{n})\|_{-1}^2.\label{projection-stability}
\end{align}
\end{lemma}

\begin{proof}
Taking $\bchi_{h}=\bhu_{j,h}^{n+1}$ in \eqref{spp-step-1}, to obtain
\begin{align}
\Bigg(\frac{\bhu_{j,h}^{n+1}-\btu_{j,h}^n}{\Delta t}, \bhu_{j,h}^{n+1}\Bigg)+\|\Bar{\nu}^\frac{1}{2}\nabla \bhu_{j,h}^{n+1}\|^2+\gamma\|\nabla\cdot\bhu_{j,h}^{n+1}\|^2+\left(2\nu_T(\hat{\bu}^{'}_{h},t^n)\nabla \bhu_{j,h}^{n+1},\nabla \bhu_{j,h}^{n+1}\right)\nonumber\\=\Big(\bif_{j}(t^{n+1}),\bhu_{j,h}^{n+1}\Big)-b^*\Big(\bhu_{j,h}^{'n}, \bhu_{j,h}^n,\bhu_{j,h}^{n+1}\Big)-\Big(\nu_j^{'}\nabla \bhu_{j,h}^{n},\nabla\bhu_{j,h}^{n+1}\Big).
\end{align}
Using polarization identity and
$(2\nu_T(\bhu^{'}_{h},t^n)\nabla \bhu_{j,h}^{n+1},\nabla \bhu_{j,h}^{n+1})=2\mu\Delta t\|\hl^{n}\nabla \bu_{j,h}^{n+1}\|^2$, we get 
\begin{align}
    \frac{1}{2\Delta t}\Big(\|\bhu_{j,h}^{n+1}\|^2-\|\btu_{j,h}^n\|^2+\|\bhu_{j,h}^{n+1}-\btu_{j,h}^n\|^2\Big)+\|\Bar{\nu}^\frac{1}{2}\nabla \bhu_{j,h}^{n+1}\|^2+\gamma\|\nabla\cdot\bhu_{j,h}^{n+1}\|^2\nonumber\\+2\mu\Delta t\|\hl^{n}\nabla \bhu_{j,h}^{n+1}\|^2=\Big(\bif_{j}(t^{n+1}),\bhu_{j,h}^{n+1}\Big)-b^*\Big(\bhu_{j,h}^{'n}, \bhu_{j,h}^n,\bhu_{j,h}^{n+1}\Big)-\Big(\nu_j^{'}\nabla \bhu_{j,h}^{n},\nabla\bhu_{j,h}^{n+1}\Big).\label{pol11}
\end{align}
Applying Cauchy-Schwarz and Young's inequalities on the forcing term, yields
  \begin{align*}
      (\bif_{j}(t^{n+1}),\bhu_{j,h}^{n+1})\le\|\bif_{j}(t^{n+1})\|_{-1}\|\nabla\bhu_{j,h}^{n+1}\|\le \frac{\alpha_j}{4}\|\nabla\bhu_{j,h}^{n+1}\|^2+\frac{1}{\alpha_j}\|\bif_{j}(t^{n+1})\|_{-1}^2.
  \end{align*}
Rewrite the above trilinear form, use \eqref{trilinear-identitiy}, Cauchy-Schwarz, H\"older's, \eqref{basic-ineq}, and Poincar\'e  inequalities, provides
\begin{align}
    -b^*\left(\bhu_{j,h}^{'n}, \bhu_{j,h}^n,\bhu_{j,h}^{n+1}\right)&=b^*\left(\bhu_{j,h}^{'n}, \bhu_{j,h}^{n+1},\bhu_{j,h}^n\right)=\left(\bhu_{j,h}^{'n}\cdot\nabla \bhu_{j,h}^{n+1}, \bhu_{j,h}^n\right)+\frac12\left(\nabla\cdot\bhu_{j,h}^{'n},\bhu_{j,h}^{n+1}\cdot \bhu_{j,h}^n\right)\nonumber\\&\le\|\bhu_{j,h}^{'n}\cdot\nabla
\bhu_{j,h}^{n+1}\|\|\bhu_{j,h}^n\|+\frac12\|\nabla\cdot\bhu_{j,h}^{'n}\|_{L^\infty}\|\bhu_{j,h}^{n+1}\|\|\bhu_{j,h}^n\|\nonumber\\&\le C\||\bhu_{j,h}^{'n}|\nabla
\bhu_{j,h}^{n+1}\|\|\nabla\bhu_{j,h}^n\|+C\|\nabla\cdot\bhu_{j,h}^{'n}\|_{L^\infty}\|\nabla\bhu_{j,h}^{n+1}\|\|\nabla\bhu_{j,h}^n\|. \label{fluc-bound-first-hat}
\end{align}
Using \eqref{eddy-viscosity}, and Young's in \eqref{fluc-bound-first-hat}, gives
  \begin{align}
      -b^*\Big(\bhu_{j,h}^{'n}, 
			\bhu_{j,h}^n,&\bhu_{j,h}^{n+1}\Big)
   \le\frac{\alpha_j}{4}\|\nabla\bhu_{j,h}^{n+1}\|^2+C\|\hl^{n}\nabla 			\bhu_{j,h}^{n+1}\|\|\nabla\bhu_{j,h}^n\|+\frac{C}{\alpha_j}\|\nabla\cdot\bhu_{j,h}^{'n}\|_{L^\infty}^2\|\nabla\bhu_{j,h}^{n}\|^2\nonumber\\&\le\frac{\alpha_j}{4}\|\nabla\bhu_{j,h}^{n+1}\|^2+\mu\Delta t\|\hl^{n}\nabla 		\bhu_{j,h}^{n+1}\|^2+\left(\frac{C}{\mu\Delta t}+\frac{C}{\alpha_j}\|\nabla\cdot\bhu_{j,h}^{'n}\|_{L^\infty}^2\right)\|\nabla\bhu_{j,h}^n\|^2.\label{fluc-bound-2}
  \end{align}
  Use of H\"older's, and Young's inequalities, gives
\begin{align*}
    -(\nu_j^{'}\nabla \bhu_{j,h}^{n},\nabla\bhu_{j,h}^{n+1})\le\|\nu_j^{'}\|_{\infty}\|\nabla \bhu_{j,h}^{n}\|\|\nabla\bhu_{j,h}^{n+1}\|\le\frac{\|\nu_j^{'}\|_{\infty}}{2}\|\nabla \bhu_{j,h}^{n}\|^2+\frac{\|\nu_j^{'}\|_{\infty}}{2}\|\nabla \bhu_{j,h}^{n+1}\|^2.
\end{align*}

 Using the above bounds, and reducing the equation \eqref{pol11}, becomes
\begin{align}
    \frac{1}{2\Delta t}\Big(\|\bhu_{j,h}^{n+1}\|^2-\|\btu_{j,h}^n\|^2+\|\bhu_{j,h}^{n+1}-\btu_{j,h}^n\|^2\Big)+\frac{\Bar{\nu}_{\min}}{2}\|\nabla \bhu_{j,h}^{n+1}\|^2+\gamma\|\nabla\cdot\bhu_{j,h}^{n+1}\|^2\nonumber\\+\mu\Delta t\|\hl^{n}\nabla \bhu_{j,h}^{n+1}\|^2\le\left(\frac{\|\nu_j^{'}\|_{\infty}}{2}+\frac{C}{\mu\Delta t}+\frac{C}{\alpha_j}\|\nabla\cdot\bhu_{j,h}^{'n}\|_{L^\infty}^2\right)\|\nabla\bhu_{j,h}^n\|^2+\frac{1}{\alpha_j}\|\bif_{j}(t^{n+1})\|_{-1}^2.
\end{align}
Dropping non-negative term from the Left-Hand-Side (LHS), and rearranging
\begin{align}
    \frac{1}{2\Delta t}\Big(\|\bhu_{j,h}^{n+1}\|^2-\|\btu_{j,h}^n\|^2\Big)+\frac{\Bar{\nu}_{\min}}{2}\left(\|\nabla \bhu_{j,h}^{n+1}\|^2-\|\nabla \bhu_{j,h}^{n}\|^2\right)\nonumber\\+\left(\frac{\alpha_j}{2}-\frac{C}{\mu\Delta t}-\frac{C}{\alpha_j}\|\nabla\cdot\bhu_{j,h}^{'n}\|_{L^\infty}^2\right)\|\nabla\bhu_{j,h}^n\|^2+\gamma\|\nabla\cdot\bhu_{j,h}^{n+1}\|^2\le \frac{1}{\alpha_j}\|\bif_{j}(t^{n+1})\|_{-1}^2.\label{estimate-b4-tild}
\end{align}
Now in \eqref{spp-step-2-1}, choose $\bv_{h}=\btu_{j,h}^{n+1}$ , and $q_{h}=\hp_{j,h}^{n+1}$ for all $n=0,1,2,\cdots,M-1$, then use of Cauchy-Schwarz, and Young’s inequalities, gives
\begin{align*}
\|\btu_{j,h}^{n+1}\|^2\le\|\bhu_{j,h}^{n+1}\|^2.
\end{align*}



Using this estimate into \eqref{estimate-b4-tild}, assuming $\alpha_j>C\|\nabla\cdot\bhu_{j,h}^{'n}\|_{L^\infty}$, choosing $\mu\ge \frac{\alpha_j}{\Delta t\left(C\alpha_j^2-\|\nabla\cdot\bhu_{j,h}^{'n}\|_{L^\infty}^2\right)}$, and dropping non-negative terms from LHS, results in 
\begin{align}
    \frac{1}{2\Delta t}\Big(\|\bhu_{j,h}^{n+1}\|^2-\|\bhu_{j,h}^n\|^2\Big)+\frac{\Bar{\nu}_{\min}}{2}\Big(\|\nabla \bhu_{j,h}^{n+1}\|^2-\|\nabla \bhu_{j,h}^{n}\|^2\Big)+\gamma\|\nabla\cdot\bhu_{j,h}^{n+1}\|^2\le \frac{1}{\alpha_j}\|\bif_{j}(t^{n+1})\|_{-1}^2.
\end{align}
Finally, multiplying both sides by $2\Delta t$, and summing over the time steps, completes the proof.
\end{proof}

\begin{remark}
   Using Lemma \ref{assump1} and triangle inequality, we can prove that for a fixed $h$, and sufficiently small $\Delta t$, $\|\nabla\cdot\bhu_{j,h}^{'n}\|_{L^\infty}\rightarrow 0$ as $\gamma\rightarrow\infty$. Numerical outcomes that support this claim are presented in Table \ref{convergence-tab-div}.
\end{remark}

\begin{remark}
    Since the SPP-EEV Algorithm \ref{SPP-FEM} is linear and finite dimensional, the stability implies the well-posedness.
\end{remark}

We want to prove the SPP-EEV Algorithm \ref{SPP-FEM} tends to the Coupled-EEV Algorithm \ref{coupled-alg} as $\gamma\rightarrow\infty$. To this end, define an orthogonal complement $\bR_h:=\bV_h^\perp\subset\bX_h$ of $\bV_h$ subject to the $\bH^1(\cD)$ norm. 

\begin{lemma}\label{CR-lemma}
If $(\bX_h,Q_h)\subset(\bX,Q)$ satisfies the \textit{LBB condition} \eqref{infsup} and  $\nabla\cdot\bX_h\subset Q_h$, then there exists a constant $C_R$ independent of $h$ such that $$\|\nabla\bv_h\|\le C_R\|\nabla\cdot\bv_h\|,\hspace{3mm}\forall\bv_h\in \bR_h.$$
\end{lemma}
\begin{proof}
See \cite{GR86, linke2017connection}
\end{proof}

\begin{lemma}\label{assump1} For a fixed $\Delta t$, and  sufficiently small $h$, if the true solution is smooth enough, then $\exists\; C_*\in\mathbb{R}^+$ (which does not dependent on $h$, and $\Delta t$), such that as $\gamma\rightarrow\infty$, 
\begin{align}
    \max_{1\le n\le M}\|\bhu_{j,h}^n\|_{L^\infty}\le C_*,\hspace{2mm}\text{and}\;\;\|\nabla\cdot\bhu_{j,h}^{n}\|_{L^\infty}\rightarrow 0, \hspace{2mm}\text{for all}\hspace{2mm}j=1,2,\cdots,J.
\end{align}
\end{lemma}
\begin{proof}
    The proof is given in Appendix \ref{appendix-B}
\end{proof}

 To analyze the convergence of the Algorithm \ref{SPP-FEM}, since $\bX_h\subset\bY_h$, we can choose $\bv_{h}=\bchi_{h}$ in \eqref{spp-step-2-1}, and combine \eqref{spp-step-2-1}-\eqref{spp-step-2-2} into \eqref{spp-step-1}, to get

\begin{align}
&\Big(\frac{\bhu_{j,h}^{n+1}-\bhu_{j,h}^n}{\Delta t}, \;\bchi_{h}\Big)+b^*\big(\hspace{-1mm}<\bhu_h>^n, \bhu_{j,h}^{n+1},\bchi_{h}\big)+\big(\Bar{\nu}\nabla \bhu_{j,h}^{n+1},\nabla \bchi_{h}\big)+\gamma\big(\nabla\cdot\bhu_{j,h}^{n+1},\nabla\cdot\bchi_{h}\big)-\big(\hp_{j,h}^{n},\nabla\cdot\bchi_{h}\big)\nonumber\\&+\Big(2\nu_T(\hat{\bu}^{'}_{h},t^n)\nabla \bhu_{j,h}^{n+1},\nabla \bchi_{h}\Big)= \big(\bif_{j}(t^{n+1}),\bchi_{h}\big)-b^*\big(\bhu_{j,h}^{'n}, \bhu_{j,h}^n,\bchi_{h}\big)-\big(\nu_j^{'}\nabla \bhu_{j,h}^{n},\nabla\bchi_{h}\big).\label{hatweak1}
\end{align}

\begin{theorem}(Convergence)\label{gamma-convergence}
Let $(\bu_{j,h}^{n+1}
,p_{j,h}^{n+1})$ and $(\bhu_{j,h}^{n+1}
,\hp_{j,h}^{n+1})$ be the solutions to the Algorithm \ref{coupled-alg}, and Algorithm \ref{SPP-FEM}, respectively, for $n=0,1,\cdots,M-1$. We then have for a given $\gamma>0$, $\mu\ge\frac12$, and time-step size $\Delta t<\min\limits_{\substack{1\le n\le M \\ 1\le j\le J}}\Big\{\frac{C\alpha_{\min}}{\|\nabla\cdot\bhu_{j,h}^{'n}\|_{L^\infty}^2}\Big\}$:
\begin{align}
    \Big(\Delta t\sum_{n=1}^M\|\nabla\hspace{-1mm}\lab\bu_h\rab^n&-\nabla\hspace{-1mm}\lab\bhu_h\rab^n\|^2\Big)^{\frac12}\le\frac{C}{\gamma} exp\lp\frac{C}{\alpha_{\min}} \left(1+\frac{\Delta t}{h^3}\right)\rp\lp  \Delta t\sum_{n=0}^{M-1}\sum_{j=1}^J\|p_{j,h}^{n+1}-\hp_{j,h}^n\|^2\rp^\frac12\nonumber\\&\times\Bigg[1+\frac{1}{\alpha_{\min}}exp\left(\frac{C}{\alpha_{\min}h^2}+\frac{C}{\Delta t}\right)\left(\frac{1}{\alpha_{\min}^2\Delta t}+\frac{1}{\Delta t}+\Delta t\right)\Bigg]^\frac12.\label{convergence-thm}
\end{align}
\end{theorem}

\begin{remark} It is obvious from \eqref{convergence-thm} that for some fixed $\Delta t$, and $h$ if $\gamma\rightarrow\infty$ the SPP-EEV scheme converges to the Coupled-EEV scheme linearly.
\end{remark}
\begin{proof} We consider the following $H^1$-orthogonal decomposition of the error $\be_{j}^{n+1}:=\bu_{j,h}^{n+1}-\bhu_{j,h}^{n+1}$ 
with $\be_{j,0}^{n+1}\in\bV_h$, and $\be_{j,\bR}^{n+1}\in\bR_h$, for $n=0,1,\cdots,M-1$.\\
\textbf{Step 1:} Estimate of $\be_{j,\bR}^{n+1}$: Subtracting the equation \eqref{hatweak1} from \eqref{couple-eqn-1}, produces to
\begin{align}
    &\frac{1}{\Delta t}\Big(\be_j^{n+1}-\be_j^n,\bchi_{h}\Big)+\left(\Bar{\nu}\nabla \be_j^{n+1},\nabla \bchi_{h}\right)+\gamma\Big(\nabla\cdot\be_{j,\bR}^{n+1},\nabla\cdot\bchi_{h}\Big)+b^*\Big(\hspace{-1.5mm}\lab\bhu_h\rab^n,\be^{n+1}_j,\bchi_{h}\Big)\nonumber\\&+b^*\Big(\hspace{-1.5mm}\lab\be\rab^n,\bu_{j,h}^{n+1},\bchi_{h}\Big)-\Big(p_{j,h}^{n+1}-\hp_{j,h}^n,\nabla\cdot\bchi_{h}\Big)+2\mu\Delta t\Big(\big\{(l^n)^2-(\hl^n)^2\big\}\nabla\bu_{j,h}^{n+1},\nabla\bchi_{h}\Big)\nonumber\\&+2\mu\Delta t\Big((\hl^n)^2\nabla\be_j^{n+1},\nabla\bchi_{h}\Big)=-b^*\Big(\bhu^{'n}_{j,h},\be^n_j,\bchi_{h}\Big)-b^*\Big(\be^{'n}_{j},\bu_{j,h}^n,\bchi_{h}\Big)-\Big(\nu_j^{'}\nabla \be_{j}^{n},\nabla\bchi_{h}\Big).\label{step-1-eqn-1}
\end{align}
Take $\bchi_{h}=\be_j^{n+1}$ in \eqref{step-1-eqn-1} which gives $b^*\Big(\hspace{-1.5mm}\lab\bhu_h\rab^n,\be^{n+1}_j,\bchi_{h}\Big)=0$, and use polarization identity, to get
\begin{align}
    \frac{1}{2\Delta t}\Big(&\|\be_j^{n+1}\|^2-\|\be_j^n\|^2+\|\be_j^{n+1}-\be_j^n\|^2\Big)+\|\Bar{\nu}^\frac12\nabla \be_j^{n+1}\|^2+\gamma\|\nabla\cdot\be_{j,\bR}^{n+1}\|^2+b^*\Big(\hspace{-1.5mm}\lab\be\rab^n,\bu_{j,h}^{n+1},\be_j^{n+1}\Big)\nonumber\\&-\Big(p_{j,h}^{n+1}-\hp_{j,h}^n,\nabla\cdot\be_{j,\bR}^{n+1}\Big)+2\mu\Delta t\|\hl^n\nabla\be_j^{n+1}\|^2+2\mu\Delta t\Big(\big\{(l^n)^2-(\hl^n)^2\big\}\nabla\bu_{j,h}^{n+1},\nabla\be_j^{n+1}\Big)\nonumber\\&=-b^*\Big(\bhu^{'n}_{j,h},\be^n_j,\be_j^{n+1}\Big)-b^*\Big(\be^{'n}_{j},\bu_{j,h}^n,\be_j^{n+1}\Big)-\nu_j^{'}\Big(\nabla \be_{j}^{n},\nabla\be_j^{n+1}\Big).\label{pol-1}
\end{align}
 Rearrange the following trilinear form, use \eqref{trilinear-identitiy}, Cauchy-Schwarz, H\"older's, Poincar\'e,  and \eqref{basic-ineq} inequalities, to get
\begin{align}
b^*\Big(&\bhu^{'n}_{j,h},\be^n_j,\be_j^{n+1}\Big)= b^*\Big(\bhu^{'n}_{j,h},\be_j^{n+1},\be_j^{n+1}-\be^n_j\Big)\nonumber\\&=\left(\bhu^{'n}_{j,h}\cdot\nabla\be_j^{n+1}, \be_j^{n+1}-\be^n_j\right)+\frac12\left(\nabla\cdot\bhu^{'n}_{j,h},\be_j^{n+1}\cdot(\be_j^{n+1}-\be_j^n)\right)\nonumber\\&\le\|\bhu^{'n}_{j,h}\cdot\nabla\be_j^{n+1}\|\|\be_j^{n+1}-\be^n_j\|+\frac12\|\nabla\cdot\bhu^{'n}_{j,h}\|_{L^\infty}\|\be_j^{n+1}\|\|\be_j^{n+1}-\be_j^n\|\nonumber\\&\le \||\bhu^{'n}_{j,h}|\nabla\be_j^{n+1}\|\|\be_j^{n+1}-\be^n_j\|+C\|\nabla\cdot\bhu^{'n}_{j,h}\|_{L^\infty}\|\nabla\be_j^{n+1}\|\|\be_j^{n+1}-\be_j^n\|\nonumber\\&\le\|\hl^n\nabla\be_j^{n+1}\|\|\be_j^{n+1}-\be^n_j\|+C\|\nabla\cdot\bhu^{'n}_{j,h}\|_{L^\infty}\|\nabla\be_j^{n+1}\|\|\be_j^{n+1}-\be_j^n\|\nonumber\\&\le \frac{\alpha_j}{8}\|\nabla\be_j^{n+1}\|^2+\Delta t\|\hl^n\nabla\be_j^{n+1}\|^2+\left(\frac{1}{4\Delta t}+\frac{C}{\alpha_j}\|\nabla\cdot\bhu_{j,h}^{'n}\|_{L^\infty}^2\right)\|\be_j^{n+1}-\be_j^n\|^2.\label{before-young-error}
\end{align}
Using H\"older's, and Young’s inequalities, we obtain
\begin{align*}
    -\Big(\nu_j^{'}\nabla \be_{j}^{n},\nabla\be_j^{n+1}\Big)&\le \frac{\|\nu_j^{'}\|_\infty}{2}\Big(\|\nabla \be_{j}^{n}\|^2+\|\nabla \be_j^{n+1}\|^2\Big).
\end{align*}
Using Cauchy-Schwarz, and Young’s inequalities, we have
\begin{align*}
    \Big(p_{j,h}^{n+1}-\hp_{j,h}^n,\nabla\cdot\be_{j,\bR}^{n+1}\Big)&\le\frac{1}{2\gamma}\|p_{j,h}^{n+1}-\hp_{j,h}^n\|^2+\frac{\gamma}{2}\|\nabla\cdot\be_{j,\bR}^{n+1}\|^2.
\end{align*}
Use the trilinear estimate in \eqref{nonlinearbound3}, Lemma \ref{lemma1}, and Young's inequality, provides
\begin{align*}
b^*\Big(\hspace{-1.5mm}\lab\be\rab^n,\bu_{j,h}^{n+1},\be_j^{n+1}\Big)&\le C \|\hspace{-1.mm}\lab\be\rab^n\hspace{-1.mm}\|\left(\|\nabla\bu_{j,h}^{n+1}\|_{L^3}+\|\bu_{j,h}^{n+1}\|_{L^\infty}\right)\|\be_j^{n+1}\|\\&\le CC_*\|\hspace{-1.mm}\lab\be\rab^n\hspace{-1.mm}\|\|\nabla\be_j^{n+1}\|\\
    &\le \frac{\alpha_j}{8}\|\nabla\be_j^{n+1}\|^2+\frac{C}{\alpha_j}\|\hspace{-1.mm}\lab\be\rab^n\hspace{-1.mm}\|^2,\\
    b^*\Big(\be^{'n}_{j},\bu_{j,h}^n,\be_j^{n+1}\Big)&\le C\|\be^{'n}_{j}\|\left(\|\nabla\bu_{j,h}^n\|_{L^3}+\|\bu_{j,h}^n\|_{L^\infty}\right)\|\be_j^{n+1}\|\\&\le CC_* \|\be^{'n}_{j}\|\|\nabla\be_j^{n+1}\|\\&\le \frac{\alpha_j}{8}\|\nabla\be_j^{n+1}\|^2+\frac{C}{\alpha_j}\|\be^{'n}_{j}\|^2.
\end{align*}
For the mixing length term, apply H\"older’s and triangle inequalities, stability estimate of Algorithm \ref{coupled-alg}, Lemma \ref{lemma1}, and \ref{assump1}, Agmon’s \cite{Robinson2016Three-Dimensional}, discrete inverse, and Young's inequalities,  to obtain
\begin{align}
    2\mu\Delta t\Big(\big\{(l^n)^2-&(\hl^n)^2\big\}\nabla\bu_{j,h}^{n+1},\nabla\be_{j}^{n+1}\Big)\le 2\mu\Delta t\|(l^n)^2-(\hl^n)^2\|_{L^\infty}\|\nabla\bu_{j,h}^{n+1}\|\|\nabla\be_j^{n+1}\|\nonumber\\
    &=2\mu\Delta t \|\sum_{i=1}^J\left(|\bu_{i,h}^{'n}|^2-|\bhu_{i,h}^{'n}|^2\right)\|_{L^\infty}\|\nabla\bu_{j,h}^{n+1}\|\|\nabla\be_j^{n+1}\|\nonumber\\
    &\le 2\mu\Delta t \sum_{i=1}^J\|(\bu_{i,h}^{'n}-\bhu_{i,h}^{'n})\cdot(\bu_{i,h}^{'n}+\bhu_{i,h}^{'n})\|_{L^\infty}\|\nabla\bu_{j,h}^{n+1}\|\|\nabla\be_j^{n+1}\|\nonumber\\&\le 2\mu\Delta t \sum_{i=1}^J\|\bu_{i,h}^{'n}-\bhu_{i,h}^{'n}\|_{L^\infty}\|\bu_{i,h}^{'n}+\bhu_{i,h}^{'n}\|_{L^\infty}\|\nabla\bu_{j,h}^{n+1}\|\|\nabla\be_j^{n+1}\|\nonumber\\
    &\le C\Delta t^{\frac{1}{2}} \sum_{i=1}^J\|\be_i^{'n}\|_{L^\infty}\left(\|\bu_{i,h}^{'n}\|_{L^\infty}+\|\bhu_{i,h}^{'n}\|_{L^\infty}\right)\|\nabla\be_j^{n+1}\|\nonumber\\
    &\le C\Delta t^{\frac{1}{2}}\sum_{i=1}^J\|\be_i^n\|_{L^\infty}\|\nabla\be_j^{n+1}\|\le C\Delta t^\frac{1}{2}h^{-\frac{3}{2}}\sum_{i=1}^J\|\be_i^n\|\|\nabla\be_j^{n+1}\|\nonumber\\
    &\le\frac{\alpha_j}{8}\|\nabla\be_j^{n+1}\|^2+\frac{C\Delta t}{\alpha_jh^3}\sum_{i=1}^J\|\be_i^n\|^2.\label{mixing-length-difference}
\end{align}
Using the above estimates in \eqref{pol-1}, and reducing, produces
\begin{align}
    \frac{1}{2\Delta t}\Big(\|\be_j^{n+1}\|^2-\|\be_j^n\|^2\Big)+\left(\frac{1}{4\Delta t}-\frac{C}{\alpha_j}\|\nabla\cdot\bhu_{j,h}^{'n}\|_{L^\infty}^2\right)\|\be_j^{n+1}-\be_j^n\|^2+\frac{\Bar{\nu}_{\min}}{2}\|\nabla \be_j^{n+1}\|^2\nonumber\\+\frac{\gamma}{2}\|\nabla\cdot\be_{j,\bR}^{n+1}\|^2+(2\mu-1)\Delta t\|\hl^n\nabla\be_j^{n+1}\|^2\le \frac{\|\nu_j^{'}\|_\infty}{2}\|\nabla \be_{j}^{n}\|^2+\frac{1}{2\gamma}\|p_{j,h}^{n+1}-\hp_{j,h}^n\|^2\nonumber\\+\frac{C}{\alpha_j}\Big(\|\hspace{-1.mm}\lab\be\rab^n\hspace{-1.mm}\|^2+\|\be^{'n}_{j}\|^2\Big)+\frac{C\Delta t}{\alpha_jh^3}\sum_{i=1}^J\|\be_i^n\|^2.\label{upperbound1}
\end{align}
Choose the tuning parameter $\mu\ge\frac{1}{2}$, and time-step size  $\Delta t<\frac{C\alpha_j}{\max\limits_{1\le n\le M}\{\|\nabla\cdot\bhu_{j,h}^{'n}\|_{L^\infty}^2\}}$ and drop non-negative terms from LHS, and rearrange
\begin{align}
    \frac{1}{2\Delta t}\Big(\|\be_j^{n+1}\|^2-\|\be_j^n\|^2\Big)+\frac{\Bar{\nu}_{\min}}{2}\left(\|\nabla \be_j^{n+1}\|^2-\|\nabla \be_j^{n}\|^2\right)+\frac{\alpha_j}{2}\|\nabla \be_j^{n}\|^2+\frac{\gamma}{2}\|\nabla\cdot\be_{j,\bR}^{n+1}\|^2\nonumber\\\le \frac{1}{2\gamma}\|p_{j,h}^{n+1}-\hp_{j,h}^n\|^2+\frac{C}{\alpha_j}\Big(\|\hspace{-1.mm}\lab\be\rab^n\hspace{-1.mm}\|^2+\|\be^{'n}_{j}\|^2\Big)+\frac{C\Delta t}{\alpha_jh^3}\sum_{i=1}^J\|\be_i^n\|^2.
\end{align}
Use triangle, and Young's inequalities, then multiply both sides by $2\Delta t$, and sum over the time steps $n=0,1,\cdots,M-1$, to obtain
\begin{align}
     \|\be_j^{M}\|^2+\alpha_j\Delta t\sum_{n=1}^M\|\nabla\be_j^n\|^2+\gamma\Delta t\sum_{n=1}^{M}\|\nabla\cdot\be_{j,\bR}^{n}\|^2\le\frac{\Delta t}{\gamma}\sum_{n=0}^{M-1}\|p_{j,h}^{n+1}-\hp_{j,h}^n\|^2\nonumber\\+\left(\frac{C}{J^2\alpha_j}\Delta t+\frac{C\Delta t^2}{\alpha_jh^3}\right)\sum_{n=1}^{M-1}\sum_{j=1}^J\|\be_j^n\|^2+\frac{C}{\alpha_j}\Delta t\sum_{n=1}^{M-1}\|\be_j^n\|^2.
\end{align} 
 Summing over $j=1,2,\cdots, J$, we have
\begin{align}
\sum_{j=1}^J\|\be_j^{M}\|^2+\alpha_{\min}\Delta t\sum_{n=1}^M\sum_{j=1}^J\|\nabla\be_j^n\|^2+\gamma\Delta t\sum_{n=1}^{M}\sum_{j=1}^J\|\nabla\cdot\be_{j,\bR}^{n}\|^2\nonumber\\\le\frac{\Delta t}{\gamma}\sum_{n=0}^{M-1}\sum_{j=1}^J\|p_{j,h}^{n+1}-\hp_{j,h}^n\|^2+\Delta t\sum_{n=1}^{M-1}\frac{C}{\alpha_{\min}}\left(1+\frac{1}{J}+\frac{J\Delta t}{h^3}\right)\sum_{j=1}^J\|\be_j^n\|^2.
\end{align} 
Applying discrete Gr\"onwall inequality given in Lemma \ref{dgl}, we get 
\begin{align}
\sum_{j=1}^J\|\be_j^{M}\|^2+\alpha_{\min}\Delta t\sum_{n=1}^M\sum_{j=1}^J\|\nabla\be_j^n\|^2+\gamma\Delta t\sum_{n=1}^{M}\sum_{j=1}^J\|\nabla\cdot\be_{j,\bR}^{n}\|^2\nonumber\\\le \frac{\Delta t}{\gamma} exp\lp \frac{CT}{\alpha_{\min}}\left(1+\frac{\Delta t}{h^3}\right)\rp\sum_{n=0}^{M-1}\sum_{j=1}^J\|p_{j,h}^{n+1}-\hp_{j,h}^n\|^2.\label{after-gronwall}
\end{align}
Using Lemma \ref{CR-lemma} with \eqref{after-gronwall} yields the following bound
\begin{align}
    \Delta t\sum_{n=1}^{M}\sum_{j=1}^J\|\nabla\be_{j,\bR}^{n}\|^2\le C_R^2\Delta t\sum_{n=1}^{M}\sum_{j=1}^J\|\nabla\cdot\be_{j,\bR}^{n}\|^2\nonumber\\\le\frac{C_R^2}{\gamma^2} exp\lp \frac{C}{\alpha_{\min}}\left(1+\frac{\Delta t}{h^3}\right)\rp\lp\Delta t\sum_{n=0}^{M-1}\sum_{j=1}^J\|p_{j,h}^{n+1}-\hp_{j,h}^n\|^2\rp.\label{step-1-bound}
\end{align}
\textbf{Step 2:} Estimate of $\be_{j,0}^{n}$: To find a bound on $\Delta t\sum\limits_{n=1}^{M}\sum\limits_{j=1}^J\|\nabla\be_{j,0}^{n}\|^2,$ take $\bchi_{h}=\be_{j,0}^{n+1}$ in \eqref{step-1-eqn-1}, which yields 
\begin{align}
    \frac{1}{\Delta t}\Big(\be_j^{n+1}-\be_j^n,\be_{j,0}^{n+1}\Big)+\|\Bar{\nu}^\frac12\nabla \be_{j,0}^{n+1}\|^2+b^*\Big(\hspace{-1.5mm}\lab\bhu_h\rab^n,\be^{n+1}_{j,\bR},\be_{j,0}^{n+1}\Big)+b^*\Big(\hspace{-1.5mm}\lab\be\rab^n,\bu_{j,h}^{n+1},\be_{j,0}^{n+1}\Big)\nonumber\\+2\mu\Delta t\Big((\hl^n)^2\nabla\be_j^{n+1},\nabla\be_{j,0}^{n+1}\Big)+2\mu\Delta t\Big(\big\{(l^n)^2-(\hl^n)^2\big\}\nabla\bu_{j,h}^{n+1},\nabla\be_{j,0}^{n+1}\Big)\nonumber\\=-b^*\Big(\bhu^{'n}_{j,h},\be^n_j,\be_{j,0}^{n+1}\Big)-b^*\Big(\be^{'n}_{j},\bu_{j,h}^n,\be_{j,0}^{n+1}\Big)-\Big(\nu_j^{'}\nabla \be_{j,0}^{n},\nabla\be_{j,0}^{n+1}\Big).\label{step-2-eqn-1}
\end{align}
Using the bound in \eqref{nonlinearbound} to the first trilinear form,  and the bound in \eqref{nonlinearbound3} to the second and fourth trilinear forms of 
		\eqref{step-2-eqn-1}, to obtain
\begin{align}
    \frac{1}{\Delta t}\Big(\be_j^{n+1}-\be_j^n,\be_{j,0}^{n+1}\Big)+\|\Bar{\nu}^\frac12\nabla \be_{j,0}^{n+1}\|^2+2\mu\Delta t\|\hl^n\nabla\be_{j,0}^{n+1}\|^2\nonumber\\\le C\|\nabla\hspace{-1mm}\lab\bhu_h\rab^n\hspace{-1mm}\|\|\nabla\be^{n+1}_{j,\bR}\|\|\nabla\be_{j,0}^{n+1}\|+C\|\hspace{-1.mm}\lab\be\rab^n\hspace{-1.mm}\|\left(\|\nabla\bu_{j,h}^{n+1}\|_{L^3}+\|\bu_{j,h}^{n+1}\|_{L^\infty}\right)\|\nabla\be_{j,0}^{n+1}\|\nonumber\\-2\mu\Delta t\Big((\hl^n)^2\nabla\be_{j,\bR}^{n+1},\nabla\be_{j,0}^{n+1}\Big)-2\mu\Delta t\Big(\big\{(l^n)^2-(\hl^n)^2\big\}\nabla\bu_{j,h}^{n+1},\nabla\be_{j,0}^{n+1}\Big)\nonumber\\-b^*\Big(\bhu^{'n}_{j,h},\be^n_j,\be_{j,0}^{n+1}\Big)+C\|\be^{'n}_{j}\|\left(\|\nabla\bu_{j,h}^n\|_{L^3}+\|\bu_{j,h}^n\|_{L^\infty}\right)\|\nabla\be_{j,0}^{n+1}\|-\Big(\nu_j^{'}\nabla \be_{j,0}^{n},\nabla\be_{j,0}^{n+1}\Big).
\end{align}
Rearrange, use \eqref{trilinear-identitiy}, Cauchy-Schwarz, H\"older's, Poincar\'e,  and \eqref{basic-ineq} inequalities, then the following trilinear form, becomes
\begin{align}    -b^*\Big(\bhu^{'n}_{j,h},\be^n_j,\be_{j,0}^{n+1}\Big)&=b^*\Big(\bhu^{'n}_{j,h},\be_{j,0}^{n+1},\be^n_j\Big)=\left(\bhu^{'n}_{j,h}\cdot\nabla\be_{j,0}^{n+1},\be^n_j\right)+\frac12\left(\nabla\cdot\bhu^{'n}_{j,h},\be_{j,0}^{n+1}\cdot\be^n_j\right)\nonumber\\&\le\|\bhu^{'n}_{j,h}\cdot\nabla\be_{j,0}^{n+1}\|\|\be^n_j\|+\frac12\|\nabla\cdot\bhu^{'n}_{j,h}\|_{L^\infty}\|\be_{j,0}^{n+1}\|\|\be_j^n\|\nonumber\\&\le \|\hl^n\nabla\be_{j,0}^{n+1}\|\|\be_j^n\|+C\|\nabla\cdot\bhu^{'n}_{j,h}\|_{L^\infty}\|\nabla\be_{j,0}^{n+1}\|\|\be_j^n\|\nonumber\\&\le\frac{\Delta t}{2}\|\hl^n\nabla\be_{j,0}^{n+1}\|+\frac{\alpha_j}{10}\|\nabla\be_{j,0}^{n+1}\|^2+\left(\frac{1}{2\Delta t}+\frac{C}{\alpha_j}\|\nabla\cdot\bhu^{'n}_{j,h}\|_{L^\infty}^2\right)\|\be_j^n\|^2.\label{before-young-e0}
		\end{align}
Using H\"older’s and Young's inequalities, yields
\begin{align*}
    -\Big(\nu_j^{'}\nabla \be_{j,0}^{n},\nabla\be_{j,0}^{n+1}\Big)\le\frac{\|\nu_j^{'}\|_{\infty}}{2}\left(\|\nabla \be_{j,0}^{n}\|^2+\|\nabla\be_{j,0}^{n+1}\|^2\right).
\end{align*}
Using the above inequality, stability estimate, and Lemma \ref{lemma1}, reducing and rearranging, to get
\begin{align}
    \frac{1}{\Delta t}\Big(\be_j^{n+1}-\be_j^n,\be_{j,0}^{n+1}\Big)+\|\Bar{\nu}^\frac12\nabla \be_{j,0}^{n+1}\|^2+\left(2\mu-1\right)\Delta t\|\hl^n\nabla\be_{j,0}^{n+1}\|^2\nonumber\\\le \frac{C}{(\Bar{\nu}_{\min}\Delta t)^\frac12}\|\nabla\be^{n+1}_{j,\bR}\|\|\nabla\be_{j,0}^{n+1}\|+CC_*\|\hspace{-1.mm}\lab\be\rab^n\hspace{-1.mm}\|\|\nabla\be_{j,0}^{n+1}\|+\left(\frac{1}{2\Delta t}+\frac{C}{\alpha_j}\|\nabla\cdot\bhu^{'n}_{j,h}\|_{L^\infty}^2\right)\|\be_j^n\|^2\nonumber\\+\frac{\alpha_j}{10}\|\nabla\be_{j,0}^{n+1}\|^2+2\mu\Delta t\Big|\Big((\hl^n)^2\nabla\be_{j,\bR}^{n+1},\nabla\be_{j,0}^{n+1}\Big)\Big|+2\mu\Delta t\Big|\Big(\big\{(l^n)^2-(\hl^n)^2\big\}\nabla\bu_{j,h}^{n+1},\nabla\be_{j,0}^{n+1}\Big)\Big|\nonumber\\+CC_*\|\be^{'n}_{j}\|\|\nabla\be_{j,0}^{n+1}\|+\frac{\|\nu_j^{'}\|_{\infty}}{2}\left(\|\nabla \be_{j,0}^{n}\|^2+\|\nabla\be_{j,0}^{n+1}\|^2\right).\label{before-time-derivative-estimate}
\end{align}
Use the polarization identity, Cauchy-Schwarz, Young's and Poincar\'e's inequalities, to evaluate the time-derivative term
\begin{align*}
    \frac{1}{\Delta t}\Big(\be_j^{n+1}-\be_j^n,\be_{j,0}^{n+1}\Big)&=\frac{1}{\Delta t}\Big(\be_j^{n+1}-\be_j^n,\be_{j}^{n+1}-\be_{j,\bR}^{n+1}\Big)\\&=\frac{1}{2\Delta t}\Big(\|\be_j^{n+1}-\be_j^n\|^2+\|\be_j^{n+1}\|^2-\|\be_j^n\|^2\Big)-\frac{1}{\Delta t}\Big(\be_j^{n+1}-\be_j^n,\be_{j,\bR}^{n+1}\Big)\\&\ge \frac{1}{2\Delta t}\Big(\|\be_j^{n+1}\|^2-\|\be_j^n\|^2\Big)-\frac{C}{\Delta t}\|\nabla\be_{j,\bR}^{n+1}\|^2.
\end{align*}
Using the above estimate, Cauchy-Schwarz's, and Young's inequalities again into \eqref{before-time-derivative-estimate}, yields
\begin{align}
    \frac{1}{2\Delta t}\Big(\|\be_j^{n+1}\|^2-\|\be_j^n\|^2\Big)+\|\Bar{\nu}^\frac12\nabla \be_{j,0}^{n+1}\|^2+\left(2\mu-1\right)\Delta t\|\hl^n\nabla\be_{j,0}^{n+1}\|^2\nonumber\\\le \left(\frac{C}{\alpha_j^2\Delta t}+\frac{C}{\Delta t}\right)\|\nabla\be^{n+1}_{j,\bR}\|^2+\frac{C}{\alpha_j}\|\hspace{-1.mm}\lab\be\rab^n\hspace{-1.mm}\|^2+2\mu\Delta t\Big|\Big((\hl^n)^2\nabla\be_{j,\bR}^{n+1},\nabla\be_{j,0}^{n+1}\Big)\Big|\nonumber\\+\left(\frac{1}{2\Delta t}+\frac{C}{\alpha_j}\|\nabla\cdot\bhu^{'n}_{j,h}\|_{L^\infty}^2\right)\|\be_j^n\|^2+\frac{2\alpha_j}{5}\|\nabla\be_{j,0}^{n+1}\|^2+\frac{C}{\alpha_j}\|\be^{'n}_{j}\|^2\nonumber\\+2\mu\Delta t\Big|\Big(\big\{(l^n)^2-(\hl^n)^2\big\}\nabla\bu_{j,h}^{n+1},\nabla\be_{j,0}^{n+1}\Big)\Big|+\frac{\|\nu_j^{'}\|_{\infty}}{2}\Big(\|\nabla \be_{j,0}^{n}\|^2+\|\nabla\be_{j,0}^{n+1}\|^2\Big).\label{before-non-linear-bounds}
\end{align}
We now use Cauchy-Schwarz, and  Young's inequalities, Lemma \ref{assump1} (which holds true for sufficiently large $\gamma$), and the stability estimate, to obtain  
\begin{align*}
    2\mu\Delta t\Big|\Big((\hl^n)^2\nabla\be_{j,\bR}^{n+1},\nabla\be_{j,0}^{n+1}\Big)\Big|&\le 2\mu\Delta t\|\hl^n\nabla\be_{j,\bR}^{n+1}\|\|\hl^n\nabla\be_{j,0}^{n+1}\|\\
    &\le2\mu^2\Delta t\|\hl^n\nabla\be_{j,\bR}^{n+1}\|^2+\frac{\Delta t}{2}\|\hl^n\nabla\be_{j,0}^{n+1}\|^2\\
    &\le2\mu^2\Delta t\|\hl^n\|_{L^\infty}^2\|\nabla\be_{j,\bR}^{n+1}\|^2+\frac{\Delta t}{2}\|\hl^n\nabla\be_{j,0}^{n+1}\|^2\\
    &\le C\Delta t\|\nabla\be_{j,\bR}^{n+1}\|^2+\frac{\Delta t}{2}\|\hl^n\nabla\be_{j,0}^{n+1}\|^2.
\end{align*}
Following the same operations as applied in \eqref{mixing-length-difference}, we get
\begin{align*}
    2\mu\Delta t\Big(\big\{(l^n)^2-(\hl^n)^2\big\}\nabla\bu_{j,h}^{n+1},\nabla\be_{j,0}^{n+1}\Big)\le\frac{\alpha_j}{10}\|\nabla\be_{j,0}^{n+1}\|^2+\frac{C\Delta t}{\alpha_jh^3}\sum_{i=1}^J\|\be_i^n\|^2.
\end{align*}
Using all these above estimates into \eqref{before-non-linear-bounds}, assuming $\mu \ge\frac12$ to drop non-negative term from LHS, using triangle and Young's inequalities and reducing, then the equation \eqref{before-non-linear-bounds} becomes
\begin{align}
    \frac{1}{2\Delta t}\Big(\|\be_j^{n+1}\|^2-\|\be_j^n\|^2\Big)+\frac{\Bar{\nu}_{\min}}{2}\|\nabla \be_{j,0}^{n+1}\|^2\le C\left(\frac{1}{\alpha_j^2\Delta t}+\frac{1}{\Delta t}+\Delta t\right)\|\nabla\be^{n+1}_{j,\bR}\|^2\nonumber\\+\Big(\frac{C}{\alpha_jJ^2}+\frac{C\Delta t}{\alpha_jh^3}\Big)\sum_{i=1}^J\|\be_i^n\|^2+\left(\frac{C}{\alpha_jJ^2}+\frac{1}{2\Delta t}+\frac{C}{\alpha_j}\|\nabla\cdot\bhu^{'n}_{j,h}\|_{L^\infty}^2\right)\|\be_j^n\|^2+\frac{\|\nu_j^{'}\|_{\infty}}{2}\|\nabla \be_{j,0}^{n}\|^2.\label{e-equ-bounded}
\end{align}
Using discrete inverse inequality, Lemma \ref{assump1}, and rearranging
\begin{align}
    \frac{1}{2\Delta t}\Big(\|\be_j^{n+1}\|^2-\|\be_j^n\|^2\Big)+\frac{\Bar{\nu}_{\min}}{2}\Big(\|\nabla \be_{j,0}^{n+1}\|^2-\|\nabla \be_{j,0}^{n}\|^2\Big)\nonumber\\+\frac{\alpha_j}{2}\|\nabla \be_{j,0}^{n}\|^2\le C\left(\frac{1}{\alpha_j^2\Delta t}+\frac{1}{\Delta t}+\Delta t\right)\|\nabla\be^{n+1}_{j,\bR}\|^2\nonumber\\+\Big(\frac{C}{\alpha_jJ^2}+\frac{C\Delta t}{\alpha_jh^3}\Big)\sum_{i=1}^J\|\be_i^n\|^2+\left(\frac{C}{\alpha_jJ^2}+\frac{1}{2\Delta t}+\frac{C}{\alpha_jh^2}\right)\|\be_j^n\|^2.
\end{align}
Multiplying both sides by $2\Delta t$, and summing over the time-step $n=0,1,\cdots,M-1$, yields
\begin{align}
    \|\be_j^M\|^2+\Bar{\nu}_{\min}\Delta t\|\nabla\be_{j,0}^M\|^2+\alpha_j\Delta t\sum_{n=1}^{M-1}\|\nabla \be_{j,0}^{n}\|^2\le C\Delta t\left(\frac{1}{\alpha_j^2\Delta t}+\frac{1}{\Delta t}+\Delta t\right)\sum_{n=1}^M\|\nabla\be^{n}_{j,\bR}\|^2\nonumber\\+C\Delta t\Big(\frac{1}{\alpha_jJ^2}+\frac{\Delta t}{\alpha_jh^3}\Big)\sum_{n=1}^{M-1}\sum_{i=1}^J\|\be_i^n\|^2+C\Delta t\left(\frac{1}{\alpha_jJ^2}+\frac{1}{\Delta t}+\frac{1}{\alpha_jh^2}\right)\sum_{n=1}^{M-1}\|\be_j^n\|^2.
\end{align}
Simplify, and again sum over $j=1,2,\cdots,J$, we have
\begin{align}
    \sum_{j=1}^J\|\be_j^M\|^2+\Delta t\sum_{n=1}^{M}\alpha_{\min}\sum_{j=1}^J\|\nabla \be_{j,0}^{n}\|^2\le \Delta t \sum_{n=1}^MC\left(\frac{1}{\alpha_{\min}^2\Delta t}+\frac{1}{\Delta t}+\Delta t\right)\sum_{j=1}^J\|\nabla\be^{n}_{j,\bR}\|^2\nonumber\\+\frac{C\Delta t}{\alpha_{\min}}\Big(1+\frac{\alpha_{\min}}{\Delta t}+\frac{1}{h^2}+\frac{J\Delta t}{h^3}\Big)\sum_{n=1}^{M-1}\sum_{i=1}^J\|\be_i^n\|^2.
\end{align}
Apply the Lemma \ref{dgl} (discrete Gr\"onwall inequality), to obtain
\begin{align}
    \sum_{j=1}^J\|\be_j^M\|^2+\alpha_{\min}\Delta t\sum_{n=1}^{M}\sum_{j=1}^J\|\nabla \be_{j,0}^{n}\|^2\nonumber\\\le C exp \left(\frac{CT}{\alpha_{\min}}\Big(1+\frac{\alpha_{\min}}{\Delta t}+\frac{1}{h^2}+\frac{\Delta t}{h^3}\Big)\right)\left(1+\frac{1}{\alpha_{\min}^2}+\Delta t^2\right)\sum_{n=1}^M\sum_{j=1}^J\|\nabla\be^{n}_{j,\bR}\|^2,\label{step-2-after-gronwell}
\end{align}
and use the estimate  \eqref{step-1-bound} in \eqref{step-2-after-gronwell}, to get
\begin{align}
    \Delta t\sum_{n=1}^{M}\sum_{j=1}^J\|\nabla \be_{j,0}^{n}\|^2\le \frac{C}{\gamma^2\alpha_{\min}} exp\lp\frac{C}{\alpha_{\min}} \left(1+\frac{\alpha_{\min}}{\Delta t}+\frac{1}{h^2}+\frac{J\Delta t}{h^3}\right)\rp \nonumber\\\times\left(1+\frac{1}{\alpha_{\min}^2}+\Delta t^2\right)\lp\sum_{n=0}^{M-1}\sum_{j=1}^J\|p_{j,h}^{n+1}-\hp_{j,h}^n\|^2\rp.
\end{align}
Using triangle and Young's inequalities
\begin{align}
    \Delta t\sum_{n=1}^{M}\|\nabla \hspace{-1.mm}\lab\be_{0}\rab^{n}\|^2\le\frac{2\Delta t}{J^2}\sum_{n=1}^{M}\sum_{j=1}^J\|\nabla \be_{j,0}^{n}\|^2\le \frac{C}{\gamma^2\alpha_{\min}} exp\lp\frac{C}{\alpha_{\min}} \left(1+\frac{\alpha_{\min}}{\Delta t}+\frac{1}{h^2}+\frac{J\Delta t}{h^3}\right)\rp\nonumber\\\times\left(1+\frac{1}{\alpha_{\min}^2}+\Delta t^2\right)\lp\sum_{n=0}^{M-1}\sum_{j=1}^J\|p_{j,h}^{n+1}-\hp_{j,h}^n\|^2\rp,\label{step-2-bound-ensemble}
\end{align}
and
\begin{align}
    \Delta t\sum_{n=1}^{M}\|\nabla\hspace{-1mm}\lab\be_{\bR}\rab^{n}\|^2\le\frac{2\Delta t}{J^2}\sum_{n=1}^{M}\sum_{j=1}^J\|\nabla\be_{j,\bR}^{n}\|^2\nonumber\\\le\frac{C}{\gamma^2} exp\lp\frac{C}{\alpha_{\min}} \left(1+\frac{\Delta t}{h^3}\right)\rp\lp  \Delta t\sum_{n=0}^{M-1}\sum_{j=1}^J\|p_{j,h}^{n+1}-\hp_{j,h}^n\|^2\rp.
\end{align}

Apply triangle and Young's inequalities on $$\|\nabla\hspace{-1mm}\lab\bu_h\rab^n-\nabla\hspace{-1mm}\lab\bhu_h\rab^n\|^2$$ to obtain the desired result.

\end{proof}

\begin{remark}
    Since  $\|\nabla\cdot\bhu_{j,h}^{n}\|_{L^\infty}\rightarrow 0$ as $\gamma\rightarrow\infty$, this relaxes the time-step restriction $$\Delta t<\min\limits_{\substack{1\le n\le M \\ 1\le j\le J}}\Big\{\frac{C\alpha_{\min}}{\|\nabla\cdot\bhu_{j,h}^{'n}\|_{L^\infty}^2}\Big\}.$$
\end{remark}
\section{SCMs}\label{scm}
 To use SCMs, we consider the stochastic space $\bW:=\bL_P^2(\Omega)$. Then, the weak form of \eqref{momentum}-\eqref{incompressibility} can be represented as: Find $\bu \in \bX \otimes \bW$ and $p \in Q \otimes \bW$ which, for almost all $t\in(0,T]$, satisfy
	\begin{align}
	\mathbb{E}[(\bu_{t},\bv)]+ \mathbb{E}[(\bu\cdot\nabla \bu,\bv)]+ \mathbb{E}[(\nu\nabla \bu,\nabla \bv)]- \mathbb{E}[(p
	,\nabla\cdot \bv)]  &  = \mathbb{E}[(\bif,\bv)],&\quad\forall \bv\in \bX \otimes \bW,\\
	\mathbb{E}[(\nabla\cdot \bu,q)]  &  =0,&\quad\forall q\in Q \otimes \bW.
	\end{align}
For input data, it is common \cite{babuvska2007stochastic,gunzburger2019evolve} in UQ to have the following assumptions $\nu(\bx,\omega)=\nu(\bx,\by(\omega))$, and $\bif(\bx,t,\omega)=\bif(
\bx,t,\by(\omega))$ for a finite dimension random vector $\by(\omega)=(y_1(\omega),y_2(\omega),\cdots,y_N(\omega))\in \bGamma\subset\mathbb{R}^N$ which is distributed according to a Joint Probability Density Function (JPDF) $\rho(\by)$ in some parameter space $\bGamma=\prod\limits_{i=1}^N\Gamma_i$ with $\mathbb{E}[\by]=\textbf{0}$, and $\mathbb{V}ar[\by]=\textbf{I}_N$.  Then, the problem is to find $\bu\in \bX \otimes \bY$ and $p \in Q \otimes \bY$ which, for almost all $t\in(0,T]$, satisfy
\begin{align}
	\int_{\bGamma}  (\bu_{t},\bv)\rho(\by)d\by &+ \int_{\bGamma} (\bu\cdot\nabla \bu,\bv)\rho(\by)d\by + \int_{\bGamma}  (\nu(\bx,\by)\nabla \bu,\nabla \bv)\rho(\by)d\by \nonumber\\&  -\int_{\bGamma}  (p ,\nabla\cdot \bv)\rho(\by)d\by 
	 = \int_{\bGamma}  (\bif,\bv)\rho(\by)d\by, \qquad\qquad\qquad\qquad\qquad\forall \bv\in \bX \otimes \bY\hspace{-1mm},\label{weak_formulation_final-1}\\
	&~~~\int_{\bGamma}  (\nabla\cdot \bu,q) \rho(\by)d\by    =0, \qquad \qquad \qquad \qquad ~~~~~~~~~~~~~~~~~~~~~~~\;\;\forall q\in Q \otimes \bY\hspace{-1mm}.\label{weak_formulation_final-2}
	\end{align}

\noindent We consider the viscosity is affinely dependent on random variables as follows:\vspace{-2mm}
 \begin{align}
\nu(\bx,\by)&=c_0(\bx)+\sum_{i=1}^Nc_i(\bx)y_i. \end{align}

Recently, for the UQ of the Quantity of Interest (QoI), $\Theta$, which can be the lift, drag, and energy, SCMs were developed. In this paper, sparse grid algorithm \cite{FTW2008} is considered as SCMs in which for a given time $t^n$ and a set of sample points $\{\by_j\}_{j=1}^{J}\subset\bGamma$, we approximate the exact solution of \eqref{momentum}-\eqref{nse-initial} by solving a discrete scheme (which can be either Coupled-EEV or SPP-EEV in this work). If the space $L_\rho^2(\bGamma)$ has a basis $\{\phi_l\}_{l=1}^{N_p}$ of dimension $N_p$, then a global polynomial interpolation is constructed as below \begin{align*}
\bu_h^{sc}(t^n\hspace{-0.5ex},\bx,\by)=\sum_{l=1}^{N_{p}}c_l(t^n\hspace{-0.5ex},\bx)\phi_l(\by),
\end{align*} to approximate the discrete solution, where $c_l(t^n,\bx)$ is the $l$-th coefficient. The sparse grid algorithm considers Leja and Clenshaw--Curtis interpolation points and their associated weights $\{w^j\}_{j=1}^J$. To provide statistical information about QoI, SCMs approximate as below $$\mathbb{E}[\Theta(\bu(t^n))]=\int_{\Gamma} \Theta(\bu(t^n),\by)\rho(\by)dy\approx\sum_{j=1}^{J}w^j\Theta(\bu_{j,h}^n).$$
 For large-scale problems with large-dimensional random inputs, SCMs are far more efficient than typical MC methods since the rate of convergence of MC in these situations results in prohibitive computational costs. In Algorithm \ref{SCM-SPP-algo}, a complete outline of the SCM-SPP-EEV is presented.
 \begin{algorithm}
	\caption{SCM-SPP-EEV}\label{SCM-SPP-algo}
	\begin{algorithmic}
 \Procedure  Sparse grid algorithm
		\State \textbf{Initialization:} Mesh, FE functions, $T$, $M$, $\{\by_{j}\}_{j=1}^{J}$, $\{w_j\}_{j=1}^{J}$\\
  \hspace{5mm}\textbf{Pre-compute:} $ \{\nu_{j}\}_{j=1}^{J}$,$ \{\bu_{j,h}^0\}_{j=1}^{J}$\\
  \hspace{5mm}\textbf{for} $n = 0, \ldots, M-1$ \textbf{do}\\
		\hspace{10mm}\textbf{for} $j = 1, \ldots, J$ \textbf{do}\\
		\hspace{15mm} To compute $\bhu^{n+1}_{j,h}$, solve \eqref{spp-step-1}-\eqref{spp-step-2-2}\\
  \hspace{15mm} Calculate $\Theta(\bhu^{n+1}_{j,h})$\\
  \hspace{10mm}\textbf{end for}\\
  \hspace{9mm} Estimate $\mathbb{E}[\Theta(\bu(t^{n+1}))] \approx \sum\limits_{j = 1}^{J} w^{j} \Theta(\bhu_{j,h}^{n+1})$\\
  \hspace{5mm}\textbf{end for}
		\EndProcedure
	\end{algorithmic}
\end{algorithm}

Instead of computing $\bhu^{n+1}_{j,h}$ by solving \eqref{spp-step-1}-\eqref{spp-step-2-2}, if we compute $\bu^{n+1}_{j,h}$ by solving \eqref{couple-eqn-1}-\eqref{couple-incompressibility} in Algorithm \ref{SCM-SPP-algo} and follow the rest of the steps, which will lead to the SCM-Coupled-EEV algorithm.
\section{Numerical Experiments}\label{numerical-experiment}
In this section, we present several numerical tests that verify the predicted convergence rates and show the performance of the scheme on some benchmark 2D problems where $\bx=(x_1,x_2)$. Pointwise-divergence free $(\mathbb{P}_2,\mathbb{P}_1^{disc})$ SV element on barycenter refined structured triangular meshes is used for the Coupled-EEV scheme (which is essentially the case $\gamma=0$), and $(\mathbb{P}_2,\mathbb{P}_1)$ TH element on regular triangular meshes is used for the SPP-EEV scheme.

\subsection{Convergence Rates Verification}

In the first experiment, we verify the theoretically found convergence rates beginning with the following analytical solution:
\[ {\bu}=\left(\begin{array}{c} \cos x_2+(1+e^t)\sin x_2 \\ \sin x_1+(1+e^t)\cos x_1 \end{array} \right)\;\;\text{and}\;\; \ p =\sin(x_1+x_2)(1+e^t),
\]
on domain $\cD=[0,1]^2$. Then, we introduce noise as $\bu_j=(1+k_j\epsilon)\bu$, and $p_j=(1+k_j\epsilon)p$, where $\epsilon$ is a perturbation parameter, $k_j:=(-1)^{j+1}4\lceil j/2\rceil/J$, $j=1,2,\cdots, J$, and $J=20$. We consider $\epsilon=0.01$, this will introduce $10\%$ noise in the initial condition, boundary condition, and the forcing functions. The forcing function $\bif_j$ is computed using the above synthetic data into \eqref{gov1}. We assume the viscosity $\nu$ is a continuous uniform random variable, and consider four random samples of size $J$ as $\nu\sim
\mathcal{U}(0.09, 0.11)$ with $\E[\nu]=0.1$, $\nu\sim
\mathcal{U}(0.009, 0.011)$ with $\E[\nu]=0.01$,  $\nu\sim
\mathcal{U}(0.0009, 0.0011)$ with $\E[\nu]=0.001$, and $\nu\sim
\mathcal{U}(0.00009, 0.00011)$ with $\E[\nu]=0.0001$. We consider the initial condition as $\bu_{j,h}^0=\bu_{j}(\bx,0)$. The Dirichlet boundary condition in the Coupled-EEV scheme and Step 1 in the SPP-EEV scheme is set as $\bu_{j,h}|_{\partial\cD}=\bu_j$, while in Step 2 in the SPP-EEV scheme, the normal velocity component is set to zero as the boundary condition.

\subsubsection{SPP-EEV scheme converges to the Coupled-EEV scheme as $\gamma\rightarrow\infty$}

We define the velocity, and pressure errors as $<\hspace{-1mm}\hat{\be}_{\bu}\hspace{-1mm}>:=<\hspace{-1mm}\bu_h\hspace{-1mm}>-<\hspace{-1mm}\bhu_{h}\hspace{-1mm}>$, and $<\hspace{-1mm}\hat{e}_p\hspace{-1mm}>:=<\hspace{-1mm}\mathscr{P}_h\hspace{-1mm}>-<\hspace{-1mm}\hp_{h,\gamma}\hspace{-1mm}>$, respectively, where 
\begin{align*}
    \mathscr{P}_{j,h}:=p_{j,h}-\frac{1}{Area(\cD)}\int_{\cD}p_{j,h} d\cD,\;\text{and}\;\;\hp_{j,h,\gamma}:=\hp_{j,h}-\frac{1}{Area(\cD)}\int_{\cD}\hp_{j,h}d\cD-\gamma\nabla\cdot\bhu_{j,h}.
\end{align*} That is, these errors are the difference between the outcomes of the coupled and projection schemes. 
We consider the simulation end time $T=1$, time-step size $\Delta t=T/10$, $\mu=1$, $\epsilon=0.01$, and $h=1/32$. Starting with $\gamma=0$, we successively increase $\gamma$ from 1e-2 by a factor of 10, record the errors in velocity and pressure, and compute the convergence rates, and finally present them in Table \ref{convergence-tab}. We observe that as $\gamma$ increases, the convergence rates asymptotically converge to 1, which is in excellent agreement with the theoretically predicted convergence rates in terms of $\gamma$ presented in \eqref{convergence-thm}. We also observe that as $\gamma\rightarrow\infty$, the divergence error tends to zero in Table \ref{convergence-tab-div}.

\begin{table}[!ht] 
\begin{center}
\small{\begin{tabular}{|c|c|c|c|c|c|c|c|c|}\hline
\multicolumn{9}{|c|}{Fixed $T=1$, $\Delta t=T/10$, $h=1/32$}\\\hline
$\hspace{-1mm}\epsilon=0.01\hspace{-1mm}$&\multicolumn{4}{c|}{$\mathbb{E}[\nu]=0.01$}&\multicolumn{4}{c|}{$\mathbb{E}[\nu]=0.001$}\\\hline
$\gamma$ & $\|\hspace{-1mm}<\hspace{-1mm}\hat{\be}_{\bu}\hspace{-1mm}>\hspace{-1mm}\|_{2,1}$ & rate &$\|\hspace{-1mm}<\hspace{-1mm}\hat{e}_p\hspace{-1mm}>\hspace{-1mm}\|_{2,0}$& rate &$\|\hspace{-1mm}<\hspace{-1mm}\hat{\be}_{\bu}\hspace{-1mm}>\hspace{-1mm}\|_{2,1}$ & rate  &$\|\hspace{-1mm}<\hspace{-1mm}\hat{e}_p\hspace{-1mm}>\hspace{-1mm}\|_{2,0}$& rate\\\hline
$0$ & 3.9912e-0& & 6.7215e-1& & 5.1823e-0 & &6.7359e-1&  \\\hline
1e-2 & 3.6882e-0 & 0.03   &6.6291e-1   &0.01& 4.6428e-0 & 0.05   &6.6413e-1   &0.01\\\hline
1e-1 & 2.7593e-0 & 0.13   &5.9839e-1   &0.04& 3.4584e-0 & 0.13   &5.9905e-1   &0.04\\\hline
1e-0 & 9.3147e-1 & 0.47   &3.1922e-1   &0.27& 1.0567e-0 & 0.51   &3.1895e-1   &0.27\\\hline
1e+1 & 1.5728e-1 & 0.77   &5.2694e-2   &0.78& 1.9040e-1 & 0.74   &5.2358e-2   &0.78\\\hline
1e+2 & 1.7306e-2 & 0.96   &5.5839e-3   &0.97&2.1293e-2 & 0.95   &5.5397e-3   &0.98\\\hline
1e+3 & 1.7479e-3 & 1.00   &6.0452e-4   &1.00& 2.1537e-3 & 1.00   &6.0788e-4   &0.96\\\hline
\end{tabular}}
\end{center}
\caption{\footnotesize SPP-EEV scheme converges to the Coupled-EEV scheme as $\gamma$ increases with $J=20$, and $\mu=1$.}\label{convergence-tab}
\end{table}

\begin{table}[!ht] 
\begin{center}
\small{\begin{tabular}{|c|c|c|c|c|c|c|c|}\hline
$\gamma$ &  0 &1e+1&1e+2&1e+3 &1e+4& 1e+5&1e+6 \\\hline
$\hspace{-1mm}\|\nabla\cdot<\hspace{-1mm}\bhu_h\hspace{-1mm}>\hspace{-1mm}\|_{\infty,0}\hspace{-1mm}$ &\hspace{-1mm}3.5450e+0\hspace{-1mm}&\hspace{-1mm}6.1481e-1&\hspace{-1mm} 9.1794e-2&\hspace{-1mm}9.5960e-3&\hspace{-1mm}9.6394e-4&\hspace{-1mm}9.6443e-6&\hspace{-1mm}9.6443e-7\\\hline
\end{tabular}}
\end{center}
\caption{\footnotesize Divergence errors in SPP-EEV method for $\mathbb{E}[\nu]=0.001$. $\|\nabla\cdot<\hspace{-1mm}\bhu_h\hspace{-1mm}>\hspace{-1mm}\|_{\infty,0}\rightarrow 0$ as $\gamma\rightarrow\infty$.}\label{convergence-tab-div}
\end{table}

\subsubsection{Spatial and temporal convergence of the SPP-EEV scheme}
Now, we define the error between the solution of SPP-EEV scheme and the exact solution as $<\hspace{-1mm}\be_{\bu}\hspace{-1mm}>:=<\hspace{-1mm}\bu^{true}\hspace{-1mm}>-<\hspace{-1mm}\bhu_h\hspace{-1mm}>$. The upper bound of this error is the same as it in \eqref{convergence-error} for large $\gamma$, which can be shown by using triangle inequality. To observe spatial convergence, we keep temporal error small enough and thus we fix a very short simulation end time $T=0.001$. We successively reduce the mesh width $h$ by a factor of 1/2, run the simulations, and record the errors, and convergence rates in Table \ref{sp-convergence-ep-0.01}. For $\gamma=$1e+6, we observe the second order convergence rates for all the three samples, which support our theoretical finding for the $(\mathbb{P}_2,\mathbb{P}_1)$ element.
\begin{table}[!ht] 
		\begin{center}
			\small\begin{tabular}{|c|c|c|c|c|c|c|c|c|}\hline
				\multicolumn{9}{|c|}{Spatial convergence (fixed $T=0.001$, $\Delta t =T/8$)}\\\hline
				$\hspace{-1mm}\epsilon=0.01\hspace{-1mm}$&\multicolumn{2}{c|}{\hspace{-1mm}$\mathbb{E}[\nu]=0.1$ }&\multicolumn{2}{c|}{\hspace{-1mm}$\mathbb{E}[\nu]=0.01$ }&\multicolumn{2}{c|}{$\mathbb{E}[\nu]=0.001$}&\multicolumn{2}{c|}{$\mathbb{E}[\nu]=0.0001$}\\\hline
				$h$& $\|\hspace{-1mm}<\hspace{-1mm}e_{\bu}\hspace{-1mm}>\hspace{-1mm}\|_{2,1}$ & rate  & $\|\hspace{-1mm}<\hspace{-1mm}e_{\bu}\hspace{-1mm}>\hspace{-1mm}\|_{2,1}$ & rate   
				 & 
				$\|\hspace{-1mm}<\hspace{-1mm}e_{\bu}\hspace{-1mm}>\hspace{-1mm}\|_{2,1}$ & rate & 
				$\|\hspace{-1mm}<\hspace{-1mm}e_{\bu}\hspace{-1mm}>\hspace{-1mm}\|_{2,1}$ & rate  
				 \\ \hline
$1/2$ & 4.4650e-4&&4.5123e-4 &&4.5128e-4& &4.5128e-4 & \\\hline
$1/4$ &1.1422e-4&1.97& 1.1568e-4 & 1.96 &1.1569e-4 & 1.96  & 1.1569e-4 & 1.96\\\hline
$1/8$ &2.8713e-5&1.99& 2.9134e-5 & 1.99 &2.9138e-5 & 1.99  & 2.9139e-5 & 1.99\\\hline
$1/16$ &7.1881e-6&2.00& 7.3380e-6 & 1.99 &7.3495e-6 & 1.99& 7.3516e-6 & 1.99\\\hline
$1/32$ &1.7983e-6&2.00& 1.8636e-6 & 1.98 &1.8938e-6 & 1.96 & 1.9133e-6 & 1.94\\\hline 

			\end{tabular}
		\end{center}
		\caption{\footnotesize Spatial errors and convergence rates of SPP-EEV scheme with $J=20$, $\mu=1$, and $\gamma=$1e+6.}\label{sp-convergence-ep-0.01} 
	\end{table}
 
 On the other hand, to observe temporal convergence, we keep fixed mesh size $h=1/64$, and simulation end time $T=1$ for $\epsilon=0.001$, and $0.01$. We run the simulations with various time-step sizes $\Delta t$ beginning with $T/1$ and successively reduce it by a factor of 1/2, record the errors, compute the convergence rates, and present them in Table \ref{tm-convergence-ep-0.01}. We observe the convergence rates approximately equal to 1. Since the linearized backward-Euler method is used in the proposed SPP-EEV scheme, the found temporal convergence rate is optimal and in excellent agreement with the theory for all the three samples.
\begin{table}[!ht] 
		\begin{center}
			\small\begin{tabular}{|c|c|c|c|c|c|c|c|c|}\hline
				\multicolumn{9}{|c|}{Temporal convergence (fixed $T=1$, $h =1/64$)}\\\hline
                \multirow{4}{*}{$\Delta t$}&\multicolumn4{c|}{$\epsilon=0.001$}&\multicolumn4{c|}{$\epsilon=0.01$}\\\cline{2-9}
&\multicolumn2{c|}{$\mathbb{E}[\nu]=0.1$}&\multicolumn2{c|}{$\mathbb{E}[\nu]=0.01$}&\multicolumn2{c|}{$\mathbb{E}[\nu]=0.1$}&\multicolumn{2}{c|}{$\mathbb{E}[\nu]=0.01$}\\\cline{2-9}
				 & $\|\hspace{-1mm}<\hspace{-1mm}e_{\bu}\hspace{-1mm}>\hspace{-1mm}\|_{2,1}$ & rate & $\|\hspace{-1mm}<\hspace{-1mm}e_{\bu}\hspace{-1mm}>\hspace{-1mm}\|_{2,1}$ & rate   
				&$\|\hspace{-1mm}<\hspace{-1mm}e_{\bu}\hspace{-1mm}>\hspace{-1mm}\|_{2,1}$  & rate & 
				$\|\hspace{-1mm}<\hspace{-1mm}e_{\bu}\hspace{-1mm}>\hspace{-1mm}\|_{2,1}$ & rate   
				 \\ \hline
$T/1$ &6.6789e-2& & 2.4162e-1& & 3.6329e-2&& 7.8186e-2& \\\hline
$T/2$ &2.7722e-2 & 1.27 &  9.8250e-2 & 1.30&1.3894e-2&1.39&3.2035e-2  &1.29 \\\hline
$T/4$ &1.2432e-2 & 1.16 & 4.3600e-2 & 1.17 &6.3101e-3&1.14& 1.5677e-2 & 1.03     \\\hline
$T/8$ &5.8667e-3 & 1.08    & 2.0466e-2 & 1.09 &3.0848e-3&1.03&  8.2576e-3 & 0.92     \\\hline
$T/16$ &2.8489e-3 & 1.04   &  9.9140e-3 & 1.05 &1.5386e-3&1.00&  4.3877e-3 & 0.91    \\\hline
$T/32$ &1.4056e-3 & 1.02   &  4.8872e-3 & 1.02 &7.7105e-4&1.00&  2.2993e-3 & 0.93     \\\hline
$T/64$ &6.9710e-4 & 1.01    &  2.4218e-3 & 1.01 &3.8646e-4&1.00&  1.1825e-3 & 0.96     \\\hline
$T/128$ &3.4749e-4 & 1.00   &  1.2056e-3 & 1.01 &1.9424e-4 & 0.99&  6.0128e-4 & 0.98     \\\hline
$T/256$ &1.7601e-4 & 0.98    & 6.0894e-4 & 0.99 &9.9057e-5 & 0.97&   3.0625e-4 & 0.97    \\\hline
			\end{tabular}
		\end{center}
		\caption{\footnotesize Temporal errors and convergence rates of SPP-EEV scheme for $\bu$ and $p$ with $J=20$, $\mu=1/2$, and $\gamma=$ 1e+6.}\label{tm-convergence-ep-0.01} 
	\end{table}

\subsection{Taylor Green-vortex (TGV) Problem \cite{taylor1937mechanism}} 

We consider the following closed form exact solution of \eqref{momentum}-\eqref{nse-initial}:
\[	{\bu}=e^{-2\nu t}\left(\hspace{-2mm}\begin{array}{c} \sin x_1\cos x_2\; \\ -\cos x_1\sin x_2\; \end{array} \hspace{-2mm}\right)\hspace{-1mm},\;\text{and}\; \ p =\frac14(\cos(2x_1)+\cos(2x_2))e^{-4\nu t},
	\]
together with the domain $\cD=[0,L]\times[0,L]$, and $\bif=\textbf{0}$. The time-dependent TGV problem shows decaying vortex as time grows. In this section, we consider the SNSE~\eqref{momentum}-\eqref{nse-initial} with a random viscosity $\nu({\bx},{\by})$, where ${\by}=(y_1,y_2,\cdots,y_5)\in\Gamma\subset\mathbb{R}^5$ is a random vector distributed according to a JPDF in a parameter space $\bGamma=\prod\limits_{i=1}^5\Gamma_i$ with $\Gamma_i:=[-\sqrt{3},\sqrt{3}]$, $\mathbb{E}[\by]=\textbf{0}$, and $\mathbb{V}ar[\by]=\textbf{I}_5$. We also consider $\mathbb{E}[\nu](\bx)=\frac{c}{1000}$ for a suitable $c>0$, $\mathbb{C}ov[\nu]({\bx},{\bx^{'}})=\frac{1}{1000^2}exp\left(-\frac{({\bx}-{\bx^{'}})^2}{l^2}\right)$, $L=\pi$ is the characteristic length, and $l$ is the correlation length. 
Then, the viscosity random field is represented by:
\begin{align}
\nu(\bx, \by_j)&=\frac{1}{1000}\psi(\bx, \by_j),
\end{align}where the Karhunen-Lo\'eve expansion as below:
\begin{align}
\psi(\bx, \by_j)=c+\left(\frac{\sqrt{\pi}l}{2}\right)^{\frac12}y_{j,1}(\omega)&+\sum_{k=1}^{q}\sqrt{\xi_k}\bigg(\sin\left(\frac{k\pi x_1}{L}\right)\sin\left(\frac{k\pi x_2}{L}\right)y_{j,2k}(\omega)\nonumber\\&+\cos\left(\frac{k\pi x_1}{L}\right)\cos\left(\frac{k\pi x_2}{L}\right)y_{j,2k+1}(\omega)\bigg), \label{eq:var-vis}
\end{align}
which truncates the infinite series to the first $q$ terms. The uncorrelated random variables $y_{j,k}$ have eigenvalues are equal to
$$\sqrt{\xi_k}=(\sqrt{\pi}l)^{\frac12}exp\left(-\frac{(k\pi l)^2}{8}\right).$$
For our test problem, we consider the correlation length $l=0.01$, $c=1$, $q=2$, $k=1,2,\cdots, 5$, $J=11$ stochastic collocation points, and $j=1,2,\cdots,J$. We consider the Clenshaw--Curtis sparse grid as the SCM which is generated by the software package TASMANIAN \cite{stoyanov2015tasmanian,doecode_6305} with 5D stochastic collocation points and their corresponding weights. An unstructured bary-centered refined triangular mesh that provides 45,087 dof is considered.

\begin{figure}
	\centering
	\subfloat[]{\includegraphics[width=0.35\textwidth,height=0.26\textwidth]{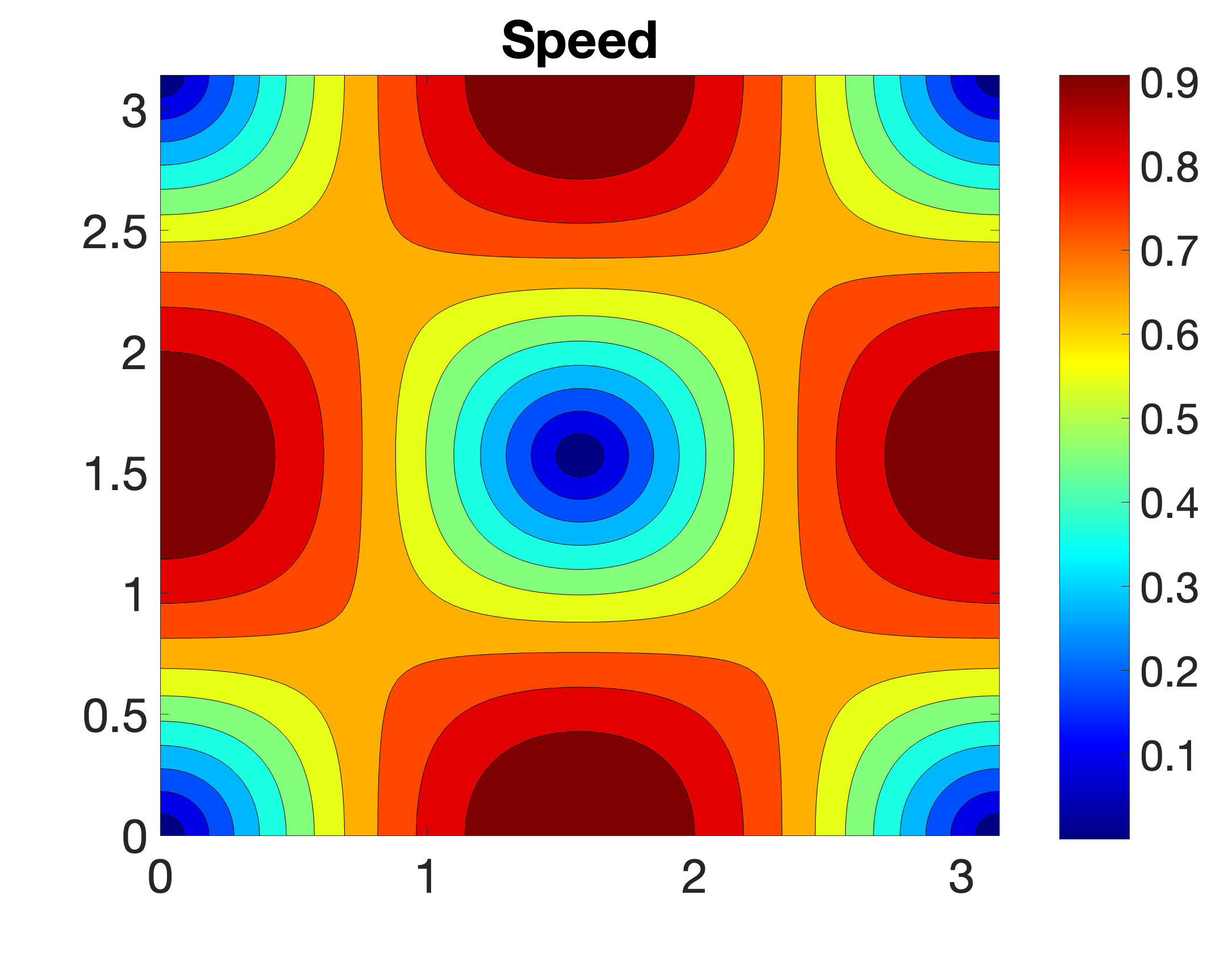}}
	\subfloat[]{\includegraphics[width=0.35\textwidth,height=0.26\textwidth]{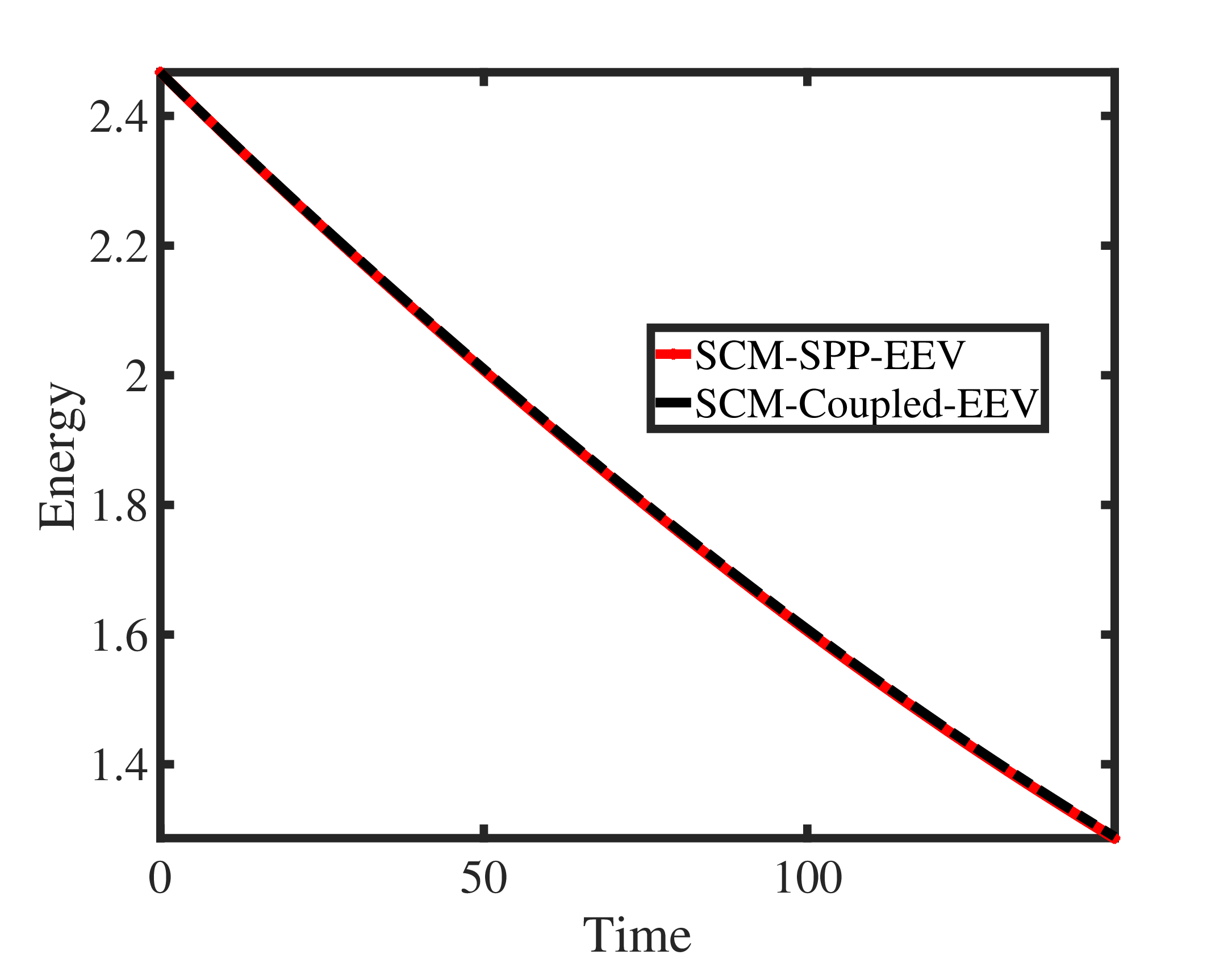}}\vspace{-4ex}
	\caption{Variable 5D random viscosity in TGV problem  for $\mathbb{E}[\nu]=0.001$: (a) Ensemble average of velocity (shown as speed contour) solution of SCM-SPP-EEV method at $t=1$, and (b) plot of Energy vs. Time for the both SCM-SPP-EEV and SCM-Coupled-EEV methods.}\vspace{-4ex}	\label{SCM-SPP-EEV-vel-pres}
\end{figure}

We consider the boundary condition $\bu_{j,h}|_{\partial\cD}=\bu$ for the SCM-Coupled-EEV and Step 1 in SCM-SPP-EEV schemes, and $\bu_{j,h}\cdot\bnh|_{\partial\cD}=0$ weakly for Step 2 in SCM-SPP-EEV scheme, and the initial condition $\bu_{j,h}^0=\bu(\bx,0)$, for both methods for $j=1,2,\cdots, J$. We run the simulations using the both methods until $T=20$ with the time-step size $\Delta t=0.1$, $\mu=1$, and $\gamma=$1e+4. We represent the approximate velocity (shown as speed) solution produced by the SCM-SPP-EEV method in Fig. \ref{SCM-SPP-EEV-vel-pres} (a) at time $t=1$.

To compare the SCM-SPP-EEV and SCM-Coupled-EEV methods, we plot their decaying Energy vs. Time graphs in Fig. \ref{SCM-SPP-EEV-vel-pres} (b) from their outcomes. For both methods, the energy at time $t=t^n$ is computed as the weighted average of $\frac12\|\hspace{-1mm}\lab\bu_h\rab^n(\bx,\by_j)\|^2$ for all stochastic collocation points. We observe an excellent agreement between the energy plots from the SCM-Coupled-EEV scheme's solution and the SCM-SPP-EEV method's solution, which supports the theory. 

\subsection{Channel Flow Over a Unit Square Step}  
\begin{figure} [ht]
		\centering	
  \subfloat[]{\includegraphics[width=0.49\textwidth,height=0.17\textwidth]{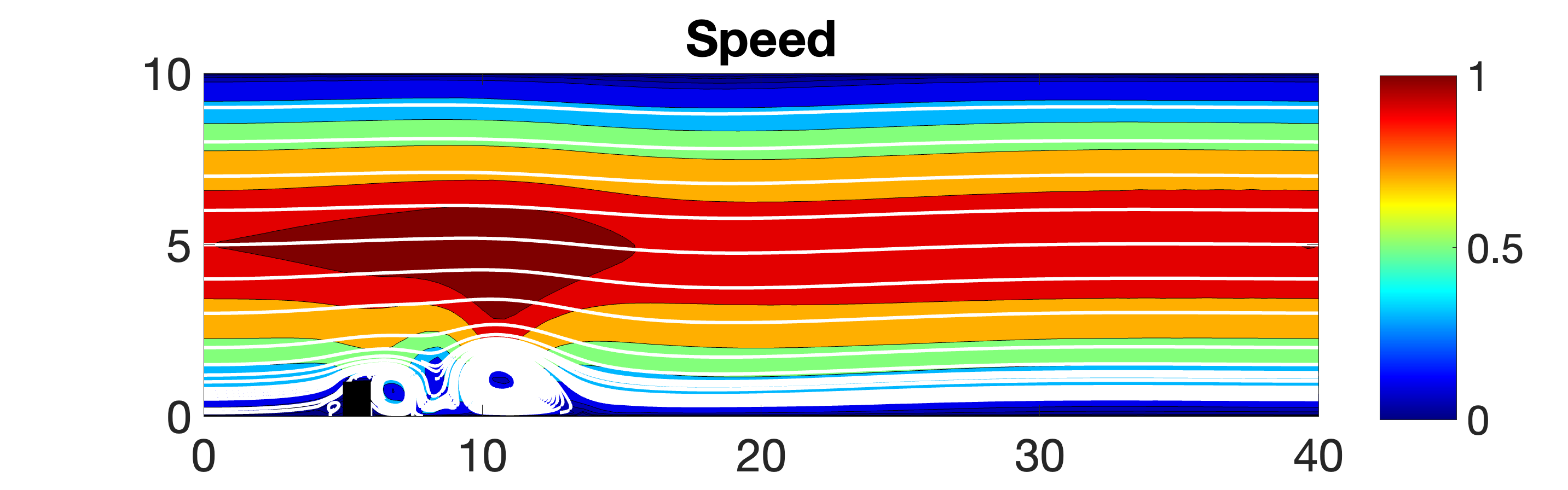}}	
  \subfloat[]{\includegraphics[width=0.3\textwidth,height=0.19\textwidth]{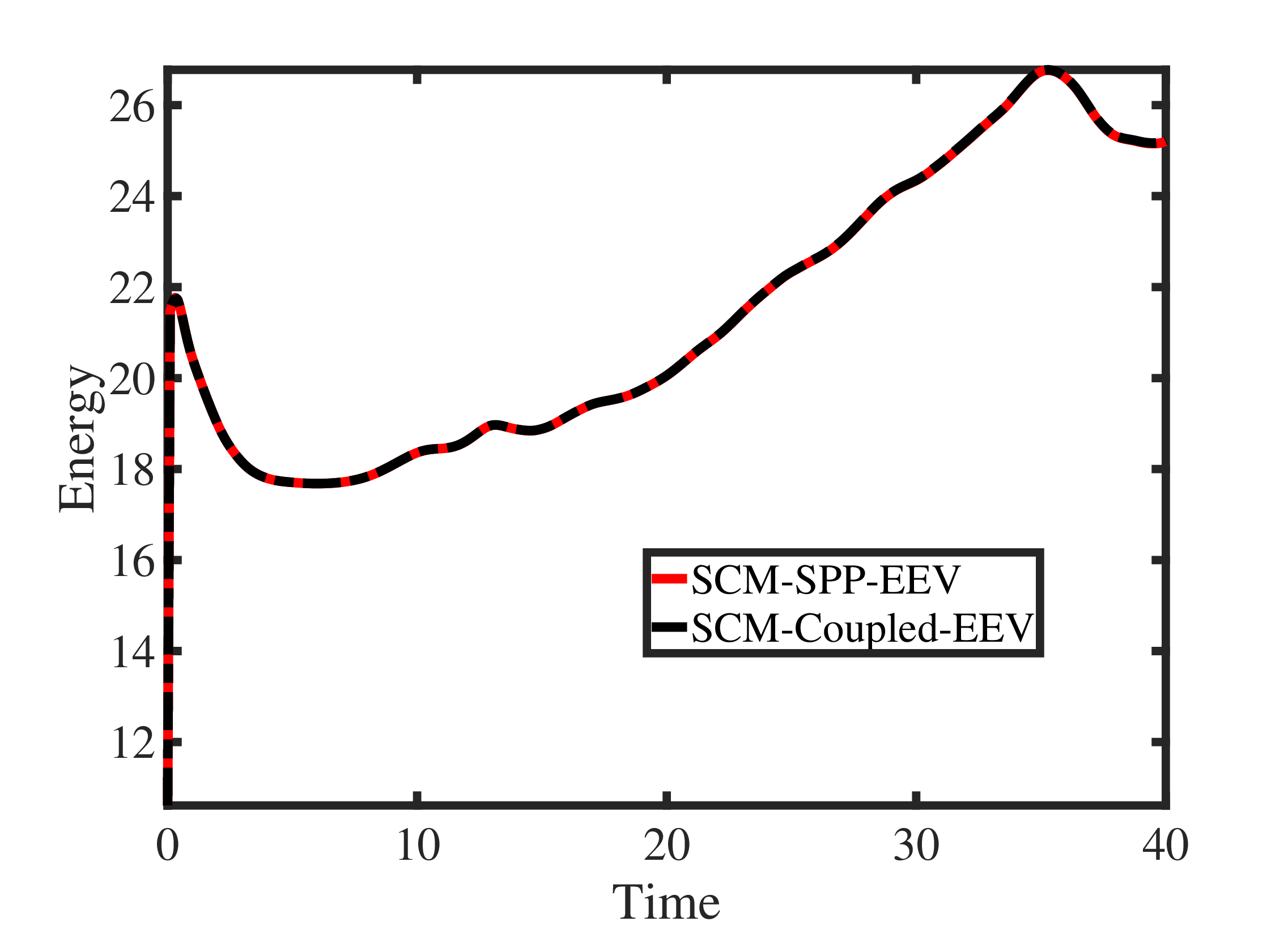}}
		\caption{\footnotesize{Variable 5D random viscosity in the flow over a step problem: (a) Ensemble average of velocity solution (shown as streamlines over the speed contour) of SCM-SPP-EEV method at $t=40$, (b) plot of Energy vs. Time for the both SCM-SPP-EEV and SCM-Coupled-EEV methods.}}\label{channel-comparison}
	\end{figure}

This experiment considers a benchmark channel flow over a unit square problem \cite{linke2017connection}. The dimension of the rectangular channel is $40\times 10$ unit$^2$, and the step is 5 units away from the inlet. The following parabolic noisy flow is considered
\begin{align*}
    \bu_{j,h}|_{inlet, outlet}=(1+k_j\epsilon){{\frac{x_2(x_2-10)}{25}}\choose{0}},
\end{align*} as inflow and outflow, where $k_j:=\frac{2j-1-J}{\lfloor \frac{J}{2}\rfloor }$, for $j=1,2,\cdots,J$, $J=11$, and $\epsilon=0.01$. No-slip boundary condition is applied to the domain walls and step for the SCM-Coupled-EEV scheme. In Step 1 of the SCM-SPP-EEV scheme, we enforce the no-slip boundary condition, and in Step 2, we set weakly, the normal velocity component vanishes on the boundary. The external force $\bif=\textbf{0}$ is considered. We start with the following initial condition
\begin{align*}
    \bu_{j,h}(\bx,0)=(1+k_j\epsilon){{\frac{x_2(x_2-10)}{25}}\choose{0}}.
\end{align*}
The random viscosity is modeled as: \begin{align}
\nu({\bx}, {\by_j})=\frac{1}{600}\psi({\bx}, {\by_j}),
\end{align}
with $\mathbb{E}[\nu](\bx)=\frac{c}{600}$, $L=40$, $y_{j}(\omega)\in[-\sqrt{3},\sqrt{3}]$, $l=0.01$, $N=5$, $c=1$, and $q=2$. The time-step size $\Delta t=0.1$, $\gamma=$ 1e+4, $\mu=1$, and run the simulation until $T=40$ using the both methods independently. The flow shows a recirculating vortex detaches behind the step \cite{layton2008numerical} as in Fig. \ref{channel-comparison} (a), which is the outcomes of the SCM-SPP-EEV method. To compare the SCM-Coupled-EEV and SCM-SPP-EEV methods, we plot the Energy vs. Time graphs in Fig. \ref{channel-comparison} (b) and found an excellent agreement between them.

\subsection{Regularized Lid-driven Cavity (RLDC) Problem}\label{RLDC}
\begin{figure}[h!] 
	\begin{center}    
          \subfloat[]
		{\includegraphics[width = 0.4\textwidth, height=0.30\textwidth]{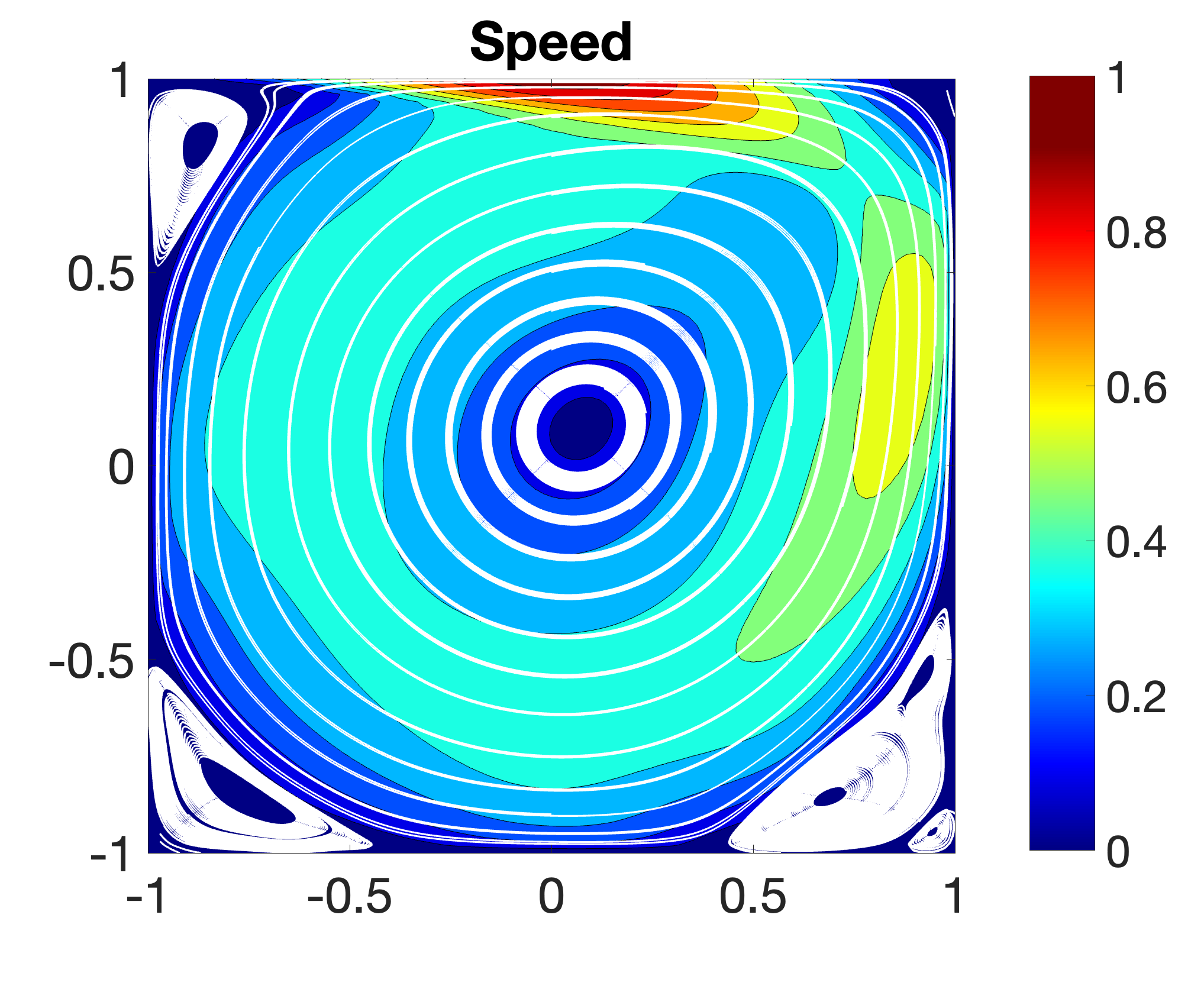}} 
  \subfloat[]{
        \includegraphics[width = 0.4\textwidth, height=0.30\textwidth]{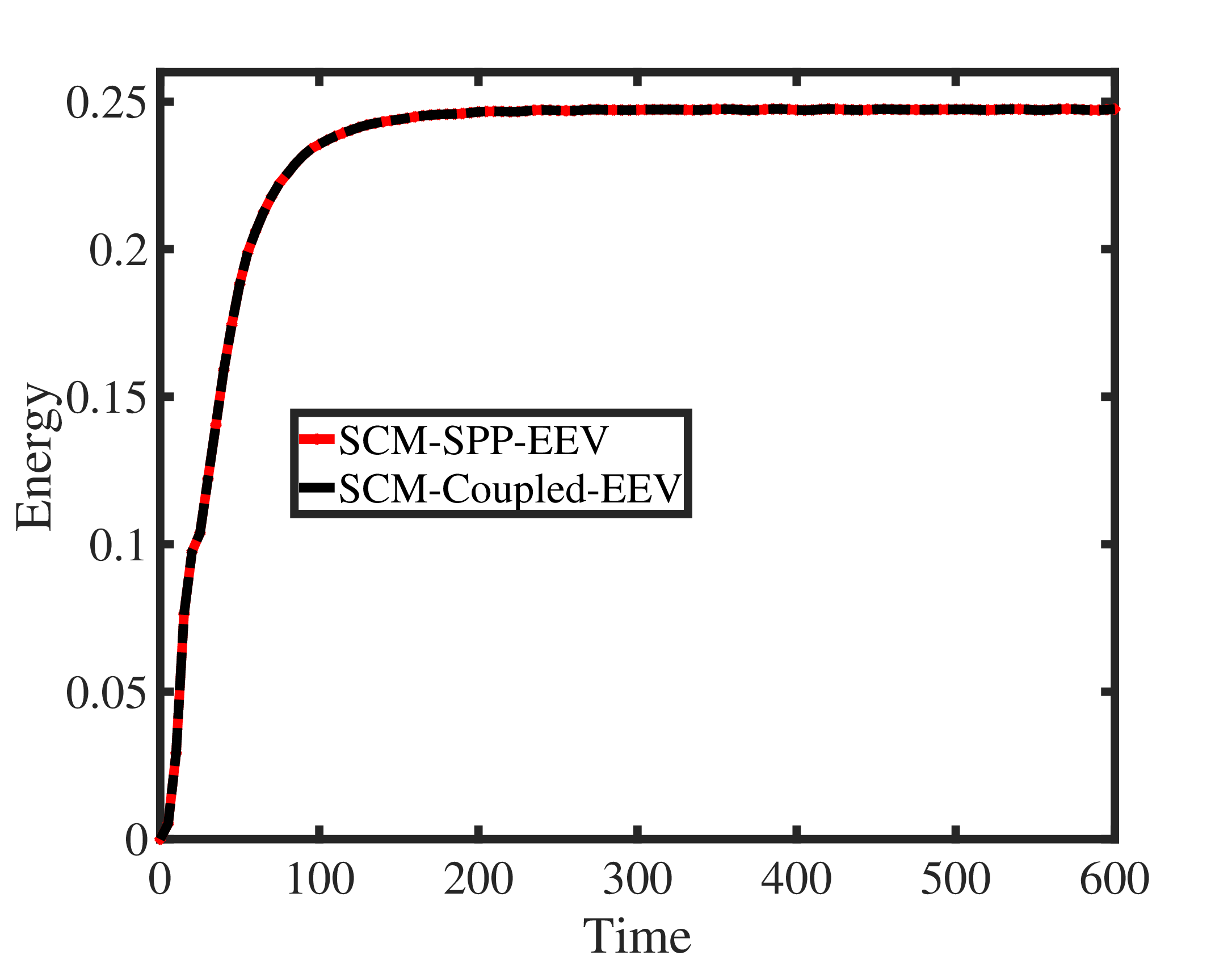}}
	\end{center}
	\caption{Variable 5D random viscosity in a RLDC problem with $\mathbb{E}[Re]=15,000$: (a) Velocity solution (shown as streamlines over the speed contour) of SCM-SPP-EEV method at $t=600$, (b) Energy vs. Time plot for both Coupled-EEV, and SCM-SPP-EEV (with $\gamma=$ 1e+4) methods.}\label{RLDC_energy_curves}
\end{figure}

A 2D benchmark RLDC problem \cite{balajewicz2013low,fick2018stabilized,lee2019study, mohebujjaman2022efficient} with a domain $\cD=(-1,1)^2$ is now under consideration. For SCM-Coupled-EEV scheme and Step 1 in SCP-SPP-EEV, all sides of the cavity are subject to no-slip boundary conditions, with the exception of the top wall (lid), where the following noise-involved boundary condition is applied:
\begin{align*}
\bu_{j,h}|_{lid}=\left (1+k_j\epsilon\right){{(1-x_1^2)^2}\choose{0}},
\end{align*}
so that the velocity of the boundary preserve the continuity. Likewise, before, we enforced the no penetration boundary condition in Step 2 in the SCM-SPP-EEV scheme. In this experiment, we model the random viscosity as below:
\begin{align}
\nu({\bx}, {\by_j})=\frac{2}{15000}\psi({\bx}, {\by_j}),
\end{align}
with $\mathbb{E}[\nu](\bx)=\frac{2c}{15000}$, $L=2$, $y_{j}(\omega)\in[-\sqrt{3},\sqrt{3}]$, $l=0.01$, $N=5$, $c=1$, and $q=2$. We conduct the simulation with an end time $T=600$ and a step size $\Delta t=5$. We consider the external force $\bif=\textbf{0}$. A perturbation parameter $\epsilon=0.01$ is applied in the boundary condition. The eddy-viscosity coefficient $\mu=1$ is considered. The unstructured triangular mesh that provides a total of 364,920 dof is considered. We represent the velocity solution of the SCM-SPP-EEV method (with $\gamma=$ 1e+4) in Fig. \ref{RLDC_energy_curves} (a) at $t=600$ and the plot of Energy vs. Time of both the SCM-Coupled-EEV, and SCM-SPP-EEV methods in Fig. \ref{RLDC_energy_curves} (b). We observe an excellent agreement of the energy plots over the time period [0, 600].

\subsection{Effect of EEV on Convection Dominated Problem}

\begin{figure}[h!] 
	\begin{center}
		\textbf{}\par\medskip            
		\includegraphics[width = 0.4\textwidth, height=0.30\textwidth]{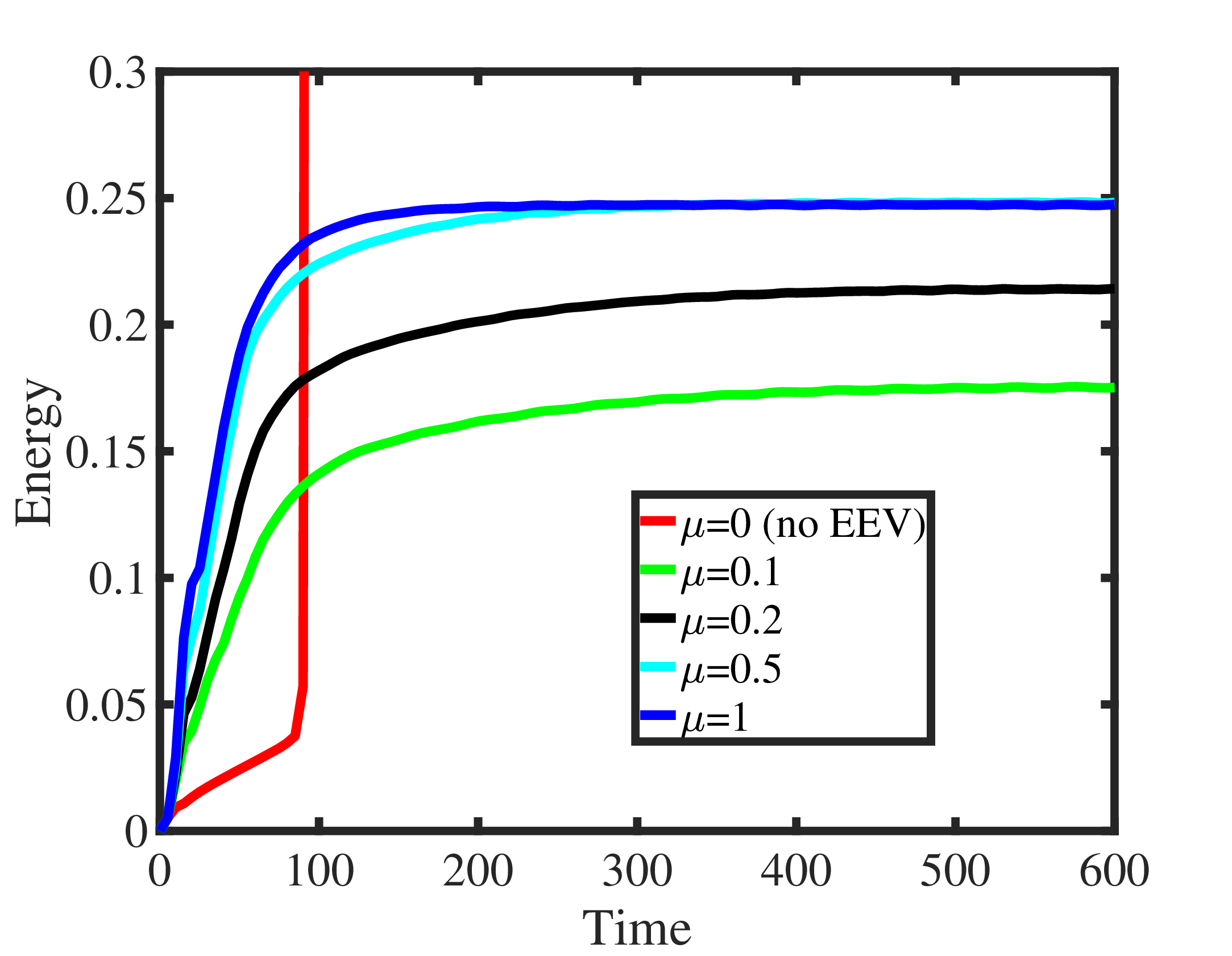}
	\end{center}
	\caption{Variable 5D random viscosity in a RLDC problem with $\mathbb{E}[Re]=15,000$: Energy vs. Time as $\mu$ varies. Solution blows up for $\mu=0$.}\label{RLDC-energy:mu-varies}
\end{figure}

The EEV based algorithms for highly ill-conditioned complex problems, e.g., RLDC problem with high $Re$, are more stable than those without it \cite{jiang2015numerical, mohebujjaman2024efficient}. To observe this, we consider the RLDC problem discussed in Section \ref{RLDC} with the same continuous and discrete model input data. We plot the Energy vs. Time graphs in Fig. \ref{RLDC-energy:mu-varies} using SCM-Coupled-EEV scheme for several values of the coefficient $\mu$ of the EEV term starting from $\mu=0$ (which is the case for without EEV algorithms). 

It is observed that the flow ensemble algorithm without EEV (setting $\mu=0$ in the SCM-Coupled-EEV scheme) blows up at around $t=90$, however, for $\mu>0$, flows remain stable until $T=600$. We also notice that as $\mu$ grows, in this case, the solution converges to the case $\mu=1$. Therefore, among the ensemble algorithms for the parameterized stochastic complex flow problems, the EEV based schemes outperform. The penalty-projection algorithm that stabilizes with EEV and grad-div terms performs better than the coupled schemes.

\section{Conclusion}
    In this work, we propose, analyze, and test the SPP-EEV fully discrete scheme which is an efficient, robust, and accurate splitting algorithm for the UQ of fluid flow problems. The unconditionally stable scheme is equipped with the following features: (i) Linear extrapolated backward-Euler method, (ii) Split by the penalty-projection method without compromising accuracy, (iii) Every subproblem is constructed in such a way that the system matrix is shared with all realizations at every time-step, (iv) To deal with convection-dominated flows, it is EEV regularized, and (v) Sparse grid SCMs is combined and proposed a more efficient and accurate SCM-SPP-EEV scheme. 

For a large penalty parameter, we demonstrated that the decoupled system converges to an equivalent linked scheme and rigorously proved the stability. The first-order temporal and optimal spatial convergence rates of the system are confirmed by a set of numerical experiments. We also verified the linear convergence of the SPP-EEV scheme to the Coupled-EEV scheme subject to the penalty parameter. Modeling the random viscosity with 5D random vectors, we showed that the SCM-SPP-EEV scheme performs well on benchmark 2D TGV, flows over a rectangular step in a channel, and RLDC problems with high expected $Re$ flow problems, and shows an excellent agreement between the SPP-EEV and SCM-Coupled-EEV algorithms in each of the cases. Finally, we showed that for convection-dominated flow, the EEV term stabilizes the flow in the RLDC problem.

In future work, following the work in \cite{linke2017connection}, we plan to propose a second-order temporal accurate efficient penalty-projection time-stepping algorithm for the UQ of NSE flow problems and analyze and test the schemes on benchmark 3D problems. Yosida-type algebraic splitting \cite{akbas2017high, rebholz2020analysis} to a block saddle point problem that arises from efficient UQ computation will be the next research avenue. 

\section{Acknowledgment} The NSF is acknowledged for supporting this research through the grant DMS-2425308. We are grateful for the generous allocation of computing time provided by the Alabama Supercomputer Authority (ASA).

\bibliographystyle{plain}
\bibliography{Penalty-Pro}

\appendix

\appendix

\section{Proof of inequality $\|\bu\cdot\nabla\bv\|\le\||\bu|\nabla \bv\|$}\label{proof-basic-eqn}
\begin{proof}
Assume $\bu=\begin{pmatrix}
    u_1\\u_2\\u_3
\end{pmatrix}$, and $\bv=\begin{pmatrix}
    v_1\\v_2\\v_3
\end{pmatrix}$, then the Euclidean norm $|\bu|=\sqrt{u_1^2+u_2^2+u_3^2}$ and
\begin{align}
\|\bu\cdot\nabla\bv\|^2&=\int_\Omega\left((\bu\cdot\nabla v_1)^2+(\bu\cdot\nabla v_2)^2+(\bu\cdot\nabla v_3)^2\right) d\Omega.\label{equn1}
\end{align}
For $\ba=\begin{pmatrix}
    a\\ b \\ c
\end{pmatrix}$, and $\bx=\begin{pmatrix}
    x\\ y \\ x
\end{pmatrix}$, using Young's inequality, we can have the following
\begin{align*}
(\ba\cdot\bx)^2=(ax+by+cz)^2&=a^2x^2+b^2y^2+c^2z^2+2abxy+2bcyz+2acxz\\&\le a^2(x^2+y^2)+b^2(x^2+y^2)+c^2z^2+2bcyz+2acxz\\&\le a^2(x^2+y^2)+b^2(x^2+y^2+z^2)+c^2(y^2+z^2)+2acxz\\&\le a^2(x^2+y^2+z^2)+b^2(x^2+y^2+z^2)+c^2(x^2+y^2+z^2)\\&=(a^2+b^2+c^2)(x^2+y^2+z^2).
\end{align*}
That is, $$(\ba\cdot\bx)^2\le |\ba|^2|\bx|^2.$$
Using the above inequality in \eqref{equn1}, we can write
\begin{align}
\|\bu\cdot\nabla\bv\|^2&\le\int_\Omega|\bu|^2\left(|\nabla v_1|^2+|\nabla v_2|^2+|\nabla v_3|^2\right)d\Omega\nonumber\\&=\int_\Omega|\bu|^2|\nabla\bv|^2d\Omega=\||\bu||\nabla\bv|\|^2=\||\bu|\nabla\bv|\|^2.\label{equn2}
\end{align}
\end{proof}

\section{Proof of Lemma \ref{lemma1}}\label{appendix-A}
Triangle inequality follows
\begin{align}
    \|\nabla \bu_{j,h}^n\|_{L^3}+\|\bu_{j,h}^n\|_{L^\infty}&\le \|\nabla(\bu_{j,h}^n-\bu_j(t^n))\|_{L^3}\nonumber\\&+\|\nabla\bu_j(t^n)\|_{L^3}+\|\bu_{j,h}^n-\bu_j(t^n)\|_{L^\infty}+\|\bu_j(t^n)\|_{L^\infty}.\label{lemma-trianlge}
\end{align}
In the RHS, 
on the first and second terms apply Sobolev embedding, and on the third, and fourth terms apply Agmon’s \cite{Robinson2016Three-Dimensional} inequality, to obtain
\begin{align}
    \|\nabla \bu_{j,h}^n\|_{L^3}+\|\bu_{j,h}^n\|_{L^\infty}&\le C\|\nabla(\bu_{j,h}^n-\bu_j(t^n))\|^{\frac12}\|\nabla^2(\bu_{j,h}^n-\bu_j(t^n))\|^{\frac12}\nonumber\\&+C\|\bu_j(t^n)\|_{H^1}^\frac12\|\bu_j(t^n)\|_{H^2}^\frac12.
\end{align}
 Apply discrete inverse inequality and regularity assumption, to get
\begin{align}
    \|\nabla \bu_{j,h}^n\|_{L^3}+\|\bu_{j,h}^n\|_{L^\infty}&\le Ch^{-\frac12}\|\nabla(\bu_{j,h}^n-\bu_j(t^n))\|+C.
\end{align}
Consider the $((P_k)^d,P_{k-1})$ element $(k\ge 2)$ for the pair $(\bu_{j,h}^n,p_{j,h}^n)$, and use the error estimate in \eqref{convergence-error}, gives 
\begin{align*}
    \|\nabla \bu_{j,h}^n\|_{L^3}+\|\bu_{j,h}^n\|_\infty&\le Ch^{-\frac12}\left(\frac{h^k}{\Delta t^{\frac12}}+\Delta t^{\frac12}+h^{1-\frac{d}{2}}\Delta t^{\frac12}+h^{k-\frac12}\right)+C.
\end{align*}
Choose $\Delta t$ so that
\begin{align*}
    \frac{h^{k-\frac12}}{\Delta t^\frac12}\le\frac{1}{C},\;
    \frac{\Delta t^\frac12}{h^\frac12}\le \frac{1}{C},\;\text{and}\;
    h^{\frac{1-d}{2}}\Delta t^\frac12\le\frac{1}{C},
\end{align*}
which gives $O(h^{2k-1})\le\Delta t\le O(h^{d-1})$, and yields
\begin{align*}
    \|\nabla \bu_{j,h}^n\|_{L^3}+\|\bu_{j,h}^n\|_\infty\le 3+C.
\end{align*}
Therefore, for $C_*:=3+C$, completes the proof.
\section{Proof of Lemma \ref{assump1}} \label{appendix-B} The proof follows the strong mathematical induction and \cite{mohebujjaman2024decoupled, wong2009analysis}.
Basic step: $\bhu_{j,h}^0=I_h(\bu_j^{true}(0,\bx)),$ where $I_h$ is an appropriate interpolation operator. Due to the regularity assumption of $\bu_j^{true}(0,\bx)$, we have $\|\bhu_{j,h}^0\|_{L^\infty}\le C_*$, for some constant $C_*>0$.\\
Inductive step: Assume for some $N\in\mathbb{N}$ and $N<M$, $\|\bhu_{j,h}^n\|_{L^\infty}\le C_*$ holds true for $n=0,1,\cdots,N$. Then, using triangle inequality and Lemma \ref{lemma1}, we have
\begin{align*}
\|\bhu_{j,h}^{N+1}\|_{L^\infty}\le\|\bhu_{j,h}^{N+1}-\bu_{j,h}^{N+1}\|_{L^\infty}+C_*.
\end{align*}
Use of Agmon’s \cite{Robinson2016Three-Dimensional} and discrete inverse inequalities, yields
\begin{align}
    \|\bhu_{j,h}^{N+1}\|_{L^\infty}\le Ch^{-\frac32}\|\bhu_{j,h}^{N+1}-\bu_{j,h}^{N+1}\|+C_*.
\end{align}
Next, using equation \eqref{after-gronwall}
\begin{align}
&\|\bhu_{j,h}^{N+1}\|_{L^\infty}\le C_*+\frac{C}{h^{\frac32}\gamma^{\frac12} } exp\lp \frac{CC_*^2}{\alpha_{\min}}+\frac{C\Delta t}{h^3\alpha_{\min}}\rp\lp  \Delta t\sum_{n=0}^{N}\sum_{j=1}^J\|p_{j,h}^{n+1}-\hp_{j,h}^n\|^2\rp^{\frac12}\hspace{-3mm}.
\end{align}

For a fixed mesh, and time-step size, as $\gamma\rightarrow \infty$, yields $\|\bhu_{j,h}^{N+1}\|_{L^\infty}\le C_*$. Hence, by the principle of strong mathematical induction, $\|\bhu_{j,h}^{n}\|_{L^\infty}\le C_*$ holds true for $0\le n\le M$.
To prove the second part, we consider the following using triangle inequality
\begin{align}
\|\nabla\cdot\bhu_{j,h}^{n}\|_{L^\infty}\le\|\nabla\cdot\left(\bu_{j,h}^{n}-\bhu_{j,h}^{n}\right)\|_{L^\infty}+\|\nabla\cdot\bu_{j,h}^{n}\|_{L^\infty}.\label{div-0}
\end{align}
Use of Agmon’s and discrete inverse inequalities, and stability estimate in \eqref{stability-couple-alg}, yields
\begin{align}
\|\nabla\cdot\bu_{j,h}^{n}\|_{L^\infty}\le Ch^{-\frac32}\|\nabla\cdot\bu_{j,h}^{n}\|\le \frac{C}{\gamma^{\frac12}\Delta t^{\frac12}h^{\frac32}}.\label{div-1}
\end{align}
Again, using Agmon’s inequality, discrete inverse inequality, and
\eqref{after-gronwall}, we have
\begin{align}
\|\nabla\cdot\left(\bu_{j,h}^{n}-\bhu_{j,h}^{n}\right)\|_{L^\infty}&\le Ch^{-\frac32}\|\nabla\cdot\left(\bu_{j,h}^{n}-\bhu_{j,h}^{n}\right)\|\le Ch^{-\frac52}\|\bu_{j,h}^{n}-\bhu_{j,h}^{n}\|\nonumber\\&\le \frac{C}{h^{\frac52}\gamma^{\frac12} } exp\lp \frac{CC_*^2}{\alpha_{\min}}+\frac{C\Delta t}{h^3\alpha_{\min}}\rp\lp  \Delta t\sum_{n=0}^{N}\sum_{j=1}^J\|p_{j,h}^{n+1}-\hp_{j,h}^n\|^2\rp^{\frac12}\hspace{-3mm}.\label{div-2}
\end{align}
Using the bounds in \eqref{div-1}-\eqref{div-2} into \eqref{div-0}, and letting $\gamma\rightarrow\infty$, complete the proof.
\end{document}